\documentclass[11pt]{article}
\usepackage{amsmath, amssymb, amsthm, mathtools}
\usepackage{bm}
\usepackage{enumerate}
\usepackage{hyperref}
\usepackage[left=1in,right=1in,top=1in,bottom=1in]{geometry}
\usepackage[backend=biber,bibencoding=utf8,maxnames=100
,doi=false,isbn=false,url=false,eprint=false,giveninits=true,date=year,sortcites]{biblatex}
\renewbibmacro{in:}{}
\usepackage{subfiles}
\usepackage{subcaption}
\usepackage{caption}
\usepackage{psfrag}
\usepackage{graphicx}
\usepackage{xkeyval}
\usepackage[noabbrev]{cleveref}
\usepackage{url}
\usepackage{xcolor}
\usepackage[per-mode=fraction]{siunitx}
\usepackage{authblk}
\usepackage{multirow}
\DeclareSIUnit\platelet{platelet}
\DeclareSIUnit\mmHg{mmHg}
\AtEveryCitekey{\clearfield{eprint}}
\date{}
\title{A Model of Fluid-Structure and Biochemical Interactions for Applications to Subclinical Leaflet Thrombosis}
\author[1,*]{Aaron Barrett}
\author[2]{Jordan A. Brown}
\author[2]{Margaret Anne Smith}
\author[3]{Andrew Woodward}
\author[4,5]{John P. Vavalle}
\author[6]{Arash Kheradvar}
\author[7,8,9,10]{Boyce E. Griffith}
\author[11]{Aaron L. Fogelson}
\affil[1]{\footnotesize Department of Mathematics, University of Utah, Salt Lake City, UT, USA}
\affil[2]{\footnotesize Department of Mathematics, University of North Carolina, Chapel Hill, NC, USA}
\affil[3]{\footnotesize Advanced Medical Imaging Lab, University of North Carolina Medical Center, Chapel Hill, NC, USA}
\affil[4]{\footnotesize University of North Carolina School of Medicine, Chapel Hill,  NC, USA}
\affil[5]{\footnotesize Division of Cardiology, Department of Medicine, University of North Carolina, Chapel Hill, NC, USA}
\affil[6]{\footnotesize Department of Biomedical Engineering, University of California Irvine, Irvine, CA, USA}
\affil[7]{\footnotesize Departments of Mathematics, Applied Physical Sciences, and Biomedical Engineering, University of North Carolina, Chapel Hill, NC, USA}
\affil[8]{\footnotesize Carolina Center for Interdisciplinary Applied Mathematics, University of North Carolina, Chapel Hill, NC, USA}
\affil[9]{\footnotesize Computational Medicine Program, University of North Carolina, Chapel Hill, NC, USA}
\affil[10]{\footnotesize McAllister Heart Institute, University of North Carolina, Chapel Hill, NC, USA}
\affil[11]{\footnotesize Departments of Mathematics and Biomedical Engineering, University of Utah, Salt Lake City, UT, USA}

\affil[*]{\footnotesize barrett@math.utah.edu}

\newcommand{\size}[1]{\left\vert #1\right\vert}

\makeatletter
\def\gsh#1{%
  \vbox{\hbox{%
    \let\\\cr
    \offinterlineskip
    \valign{&\hb@xt@2\p@{\hss$##$\hss}\vskip.2ex\cr#1\crcr}%
  }\vskip-.36ex}%
}
\def\gshsym{\@ifstar\gsh@ssym\gsh@sym}
\def\gsh@sym#1#2{\mathrlap{\overset{#1}{\phantom{#2}}}#2}
\def\gsh@ssym#1#2{\overset{#1}{#2}{\vphantom{#2}}}
\makeatother
\DeclarePairedDelimiter\norm{\lVert}{\rVert}%

\newcommand{\VV}{\mathbf{V}}

\newcommand{\uu}{\mathbf{u}}
\newcommand{\grad}{\nabla}
\newcommand{\Lap}{\Delta}
\newcommand{\ff}{\mathbf{f}}
\newcommand{\parens}[1]{\mathopen{}\left(#1\right)\mathclose{}}
\newcommand{\xx}{\mathbf{x}}

\renewcommand{\aa}{\mathbf{a}}
\newcommand{\XX}{\mathbf{X}}

\newcommand{\FF}{\mathbf{F}}

\newcommand{\bchi}{\bm{\chi}}

\newcommand{\nn}{\mathbf{n}}

\newcommand{\cc}{\mathbf{c}}

\newcommand{\etal}{et al.}
\renewcommand{\AA}{\mathbf{A}}
\newcommand{\PP}{\mathbb{P}}
\newcommand{\FFbb}{\mathbb{F}}

\newcommand{\cf}{c_\text{f}}
\newcommand{\Cb}{C_\text{b}}
\newcommand{\kon}{k_\text{on}}
\newcommand{\koff}{k_\text{off}}
\newcommand{\Cbmax}{C_\text{b}^\text{max}}
\newcommand{\Interior}{\Omega_t^{\text{f}^-}}
\newcommand{\InteriorBdry}{\Gamma_t^{\text{f}^-}}

\newcommand{\Exterior}{\Omega_t^{\text{f}^+}}
\newcommand{\Fluid}{\Omega_t^\text{f}}
\newcommand{\Solid}{\Omega_t^\text{s}}
\newcommand{\SolidRef}{\Omega_0^\text{s}}
\newcommand{\QAO}{Q_\text{Ao}}
\newcommand{\RC}{R_\text{c}}
\newcommand{\PAO}{P_\text{Ao}}
\newcommand{\PWK}{P_\text{Wk}}
\newcommand{\RP}{R_\text{p}}
\newcommand{\QLVOT}{Q_\text{LVOT}}
\newcommand{\PLVOT}{P_\text{LVOT}}
\newcommand{\RLVOT}{R_\text{LVOT}}
\newcommand{\PLV}{P_\text{LV}}
\newcommand{\CLV}{C_\text{LV}}
\newcommand{\QMV}{Q_\text{MV}}
\newcommand{\RMV}{R_\text{MV}}
\newcommand{\CLA}{C_\text{LA}}
\newcommand{\PLA}{P_\text{LA}}
\newcommand{\QVEIN}{Q_\text{vein}}
\newcommand{\Js}{J_\text{s}}

\makeatletter
\newlength{\sfp@hseplen}\newlength{\sfp@vseplen}
\define@cmdkey{subfigpos}[sfp@]{pos}[ul]{}
\define@cmdkey{subfigpos}[sfp@]{font}[\small]{}
\define@cmdkey{subfigpos}[sfp@]{vsep}[2\baselineskip]{\setlength{\sfp@vseplen}{\sfp@vsep}}
\define@cmdkey{subfigpos}[sfp@]{hsep}[10pt]{\setlength{\sfp@hseplen}{\sfp@hsep}}
\newcommand{\subfigimg}[3][,]{%
  \setkeys{Gin,subfigpos}{pos,font,vsep,hsep,#1}
  \setbox1=\hbox{\includegraphics{#3}}
  \ifnum\pdfstrcmp{\sfp@pos}{ul}=0
    \leavevmode\rlap{\usebox1}
    \rlap{\hspace*{\sfp@hsep}\raisebox{\dimexpr\ht1-\sfp@vsep}{\sfp@font{#2}}}
    \phantom{\usebox1}
  \else\ifnum\pdfstrcmp{\sfp@pos}{ur}=0
    \leavevmode\usebox1
    \llap{\raisebox{\dimexpr\ht1-\sfp@vsep}{\sfp@font{#2}}\hspace*{\sfp@hsep}}
  \else\ifnum\pdfstrcmp{\sfp@pos}{lr}=0
    \leavevmode\usebox1
    \llap{\raisebox{\sfp@vsep}{\sfp@font{#2}}\hspace*{\sfp@hsep}}
  \else
    \leavevmode\rlap{\usebox1}
    \rlap{\hspace*{\sfp@hseplen}\raisebox{\sfp@vsep}{\sfp@font{#2}}}
    \phantom{\usebox1}
  \fi\fi\fi
}
\crefname{figure}{Figure}{Figures}

\makeatother
\begin{document}
\maketitle

\begin{abstract}
Subclinical leaflet thrombosis (SLT) is a potentially serious complication of aortic valve replacement with a bioprosthetic valve in which blood clots form on the replacement valve. SLT is associated with increased risk of transient ischemic attacks and strokes and can progress to clinical leaflet thrombosis. SLT following aortic valve replacement also may be related to subsequent structural valve deterioration, which can impair the durability of the valve replacement. Because of the difficulty in clinical imaging of SLT, models are needed to determine the mechanisms of SLT and could eventually predict which patients will develop SLT. To this end, we develop methods to simulate leaflet thrombosis that combine fluid-structure interaction and a simplified thrombosis model that allows for deposition along the moving leaflets. Additionally, this model can be adapted to model deposition or absorption along other moving boundaries. We present convergence results and quantify the model's ability to realize changes in valve opening and pressures. These new approaches are an important advancement in our tools for modeling thrombosis in which they incorporate both adhesion to the surface of the moving leaflets and feedback to the fluid-structure interaction.
\end{abstract}
\section{Introduction}
Subclinical leaflet thrombosis is a potentially serious complication of bioprosthetic aortic valve replacement and may occur following either surgical or transcatheter aortic valve replacement. Although bioprosthetic heart valves (BHVs) are remarkably less thrombogenic than mechanical heart valves (MHVs), clinical valve thrombosis can occur as a life-threatening complication. Recent studies \cite{Makkar2017,Makkar2015,Schymik2015} have suggested that the rate of subclinical leaflet thrombosis (SLT) is as high as 13--38\% \cite{Rosseel2019}. SLT is associated with increased risk of transient ischemic attacks and strokes, acute myocardial infarction, and accelerated valve deterioration \cite{Ramana2019}. Further, if left untreated, SLT can progress to clinical valve thrombosis. While a cardiac computed tomography (CT) scan can detect SLT, predicting which patients will develop SLT is currently not possible. Accordingly, there is a need for of computational tools to model the fluid-structure and biochemical interactions that predispose a particular patient to develop SLT.

Prior work to model leaflet thrombosis has focused on computational fluid dynamics (CFD) simulations of blood flow through the valve. Plitman Mayo \etal\ \cite{PlitmanMayo2020} performed CFD experiments of deployed transcatheter aortic valve replacements (TAVRs) to determine areas of stagnated blood flow, suggesting possible sites of thrombosis formation. Vahidkhah \etal\ \cite{Vahidkhah2017} compared blood residence times behind the coronary and non-coronary leaflets after a TAVR procedure and determined similar residence times for all the leaflets. Kivi \etal\ \cite{Kivi2020} performed two dimensional fluid-structure interaction (FSI) simulations with leaflets of varying stiffness. A common finding in CFD and FSI simulations is the presence of stagnant regions in the aortic sinus, in which blood clots are thought to form. Hatoum \etal\ \cite{Hatoum2021} combined a CFD model of flow through patient specific geometry post-TAVR with a reduced order model that predicted thrombus growth based on the wall shear stress and percent stasis volume measurements. While they were able to determine a correlation between circulation and amount of thrombosis, they concluded that finer flow metrics or FSI analysis are needed to fully predict thrombosis.

Mathematical and computational models of thrombosis have also been developed, but methods suitable for modeling thrombosis on dynamic flexible structures, which is critical for describing leaflet thrombosis, are lacking. Fogelson \etal\ \cite{Fogelson2012,Link2018} developed a model of intravascular platelet deposition and determined the sensitivity of thrombus formation due to various chemical and platelet factors. Du and Fogelson \cite{Du2018} developed a multiphase model of platelet aggregation in which the thrombus is modeled as a viscoelastic fluid. This model can be seen as an extension of models by Fogelson and Guy \cite{Fogelson2008} that were created to study thrombus formation in a moving fluid. Models describing flowing platelets and platelet deposition onto a stationary vessel wall have been developed using a variety of multiscale modeling and computational approaches \cite{Diamond2013,Wu2014,Zhang2014,Tosenberger2016}. These models describe both fluid-phase transport of platelets and the influence of platelet deposits on the hemodynamics through and near the deposits. In these models, the platelets deposit over stationary surfaces. However, to our knowledge, no thrombosis model has yet been developed that allows for thrombus growth on a surface whose motion is determined by solving an FSI problem, e.g., a heart valve leaflet. 

There are several models that couple the advection and diffusion of chemical species and their sources from immersed boundaries \cite{Santiago2022a,Santiago2022b,Hopkins2002,}. Typically, these models use sources that are then spread from the immersed boundary to the surrounding fluid using the regularized delta function. Restricting species from diffusing across the interface remains a challenge. While many different methods have been proposed to restrict diffusion and enforce Robin boundary conditions across a moving interface \cite{Towers2018,Huang2009,Xu2006,Thirumalaisamy2022,Helgadottir2015}, there are far fewer that have tested the method in the context of an immersed boundary model. Chen and Lai \cite{Chen2014} used a diffuse domain approach to model absorption of surfactants on a surface. Their approach is based on the methods introduced by Li \etal\ \cite{Li2009b} who demonstrated that this method enforces the boundary condition at first order accuracy. In the methods used herein, we enforce a boundary condition without smoothing the interface, leading to second order accuracy up to and including the boundary \cite{Barrett2022}.

The present study introduces new numerical methods to simulate the deposition of material onto thin moving leaflets. The leaflet and fluid motion are determined through an FSI calculation and the material deposition feeds back onto the FSI calculation by modifying the leaflet's mechanical properties. While we refer to the fluid as blood and the deposited material as platelets, the current paper deals only with a prototype of a situation that would arise in modeling leaflet thrombosis. In a complete model of leaflet thrombosis, the deposited material would consist of platelets, fibrin, and potentially other inflammatory cells. The model would include mechanisms for activating platelets through contact with the leaflet surface, exposure to high shear stress, or encounter of soluble platelet agonists \cite{Fogelson2015,Du2018,Du2020}. It would also include treatment of coagulation biochemistry \cite{Link2020,Link2018,Leiderman2011,Kuharsky2001} coupled with fibrin polymerization \cite{Fogelson2022,Fogelson2010}. The current work is a major step towards simulating the dynamics of such a model.

\section{Continuous Equations}\label{sec:CEM}
We consider an FSI model of thrombus formation on the aortic valve leaflets. The valve geometry is created via slicing a three-dimensional reconstruction of a typical trileaflet aortic valve, as will be discussed further in \cref{sec:mesh}. In this simplified model, fluid phase platelets can bind to the leaflet surface while the surface-bound platelets stiffen the leaflets and can also dissociate back into the fluid.

\subsection{Fluid-Structure Interaction}
\begin{figure}
\begin{center}
\psfragscanon
\psfrag{Interior}{\scalebox{0.85}{$\Interior$}}
\psfrag{Exterior}{\scalebox{0.85}{$\Exterior$}}
\psfrag{Solid}{\scalebox{0.85}{$\Solid = \Omega\setminus\Fluid$}}
\psfrag{Gamma}{\scalebox{0.85}{$\InteriorBdry$}}
\psfrag{QQ}{\scalebox{0.85}{$\PAO$}}
\psfrag{Qao}{\scalebox{0.85}{$\QAO$}}
\psfrag{Rc}{\scalebox{0.85}{$\RC$}}
\psfrag{C}{\scalebox{0.85}{$C$}}
\psfrag{Pwk}{\scalebox{0.85}{$\PWK$}}
\psfrag{Rp}{\scalebox{0.85}{$\RP$}}
\psfrag{Plvot}{\scalebox{0.85}{$\PLVOT$}}
\psfrag{Rlvot}{\scalebox{0.85}{$\RLVOT$}}
\psfrag{Qlvot}{\scalebox{0.85}{$\QLVOT$}}
\psfrag{Qmv}{\scalebox{0.85}{$\QMV$}}
\psfrag{Rmv}{\scalebox{0.85}{$\RMV$}}
\psfrag{Clv}{\scalebox{0.85}{$\CLV$}}
\psfrag{Cla}{\scalebox{0.85}{$\CLA$}}
\psfrag{Qvein}{\scalebox{0.85}{$\QVEIN$}}
\psfrag{Plv}{\scalebox{0.85}{$\PLV$}}
\psfrag{Pla}{\scalebox{0.85}{~~$\PLA$}}
\includegraphics[width = 0.7\linewidth]{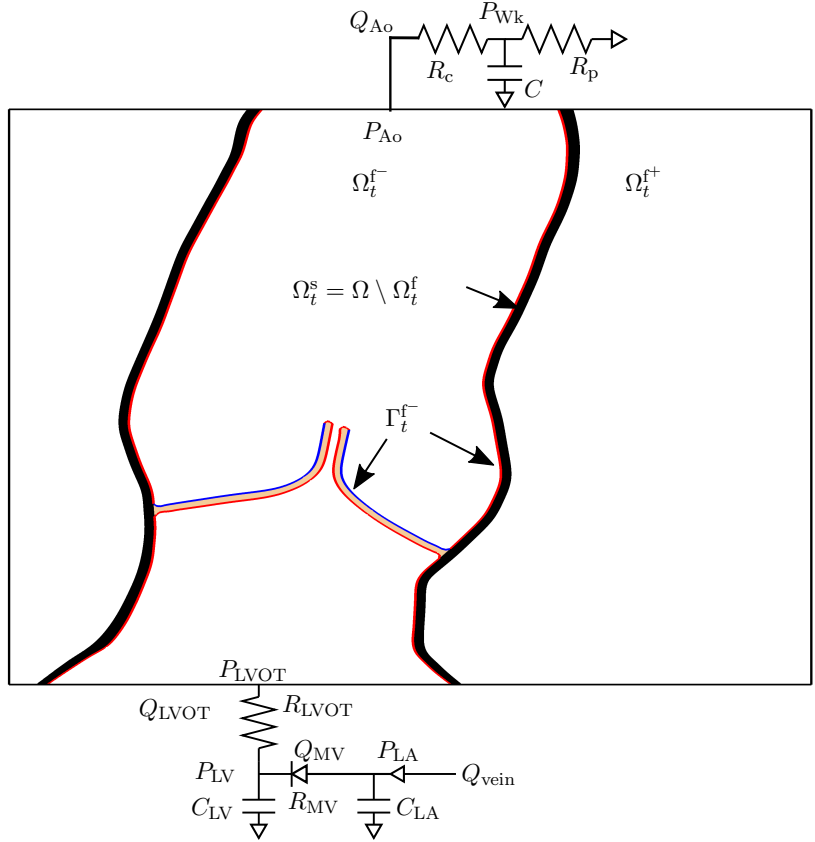}
\end{center}
\caption{The domain is decomposed into a solid subdomain $\Solid$ (denoted in the black and tan curves) and a fluid subdomain $\Fluid$, which itself is partitioned into interior subregion $\Interior$, which corresponds to the lumen, and an exterior subregion $\Exterior$, which corresponds to the space outside the vessel. The lumen boundary $\InteriorBdry = \partial\Interior$ is composed of regions where platelet binding can occur (shown in blue) and no penetration conditions for the platelets (shown in red). The inlet to the domain is the left ventricle outflow tract (LVOT) and the outlet is the ascending aorta. Boundary conditions at the inlet are determined using a time-dependent elastance based model of the heart, including the left ventricle (LV), mitral valve (MV), and the left atrium (LA). The outlet boundary conditions are determined using a three element Windkessel model \cite{Lee2021a,Griffith2009b}.}\label{fig:domains}
\end{figure}

The fluid-structure system is modeled using the immersed finite element/finite difference method \cite{Griffith2017}. In this approach, a fixed computational domain $\Omega$ is partitioned into a time-dependent fluid subdomain $\Fluid$ and a time-dependent solid subdomain $\Solid$, so that $\Omega = \Fluid \cup \Solid$. The fluid domain is further subdivided into the lumen $\Interior$ (i.e., the space occupied by the blood, in which platelets are free to advect and diffuse) and the space outside the aortic root $\Exterior$, with $\Fluid = \Interior\cap\Exterior$; see \cref{fig:domains}. We denote Eulerian physical coordinates with $\xx$. The solid domain is tracked using Lagrangian material coordinates $\XX$, and the mapping between the reference and current coordinates is $\bchi\parens{\XX,t}$. The motion of the fluid-structure system is described by 
\begin{subequations}\label{eq:fl_sl_motion}
\begin{align}
\rho\parens{\frac{\partial \uu\parens{\xx,t}}{\partial t} + \uu\parens{\xx,t}\cdot\grad\uu\parens{\xx,t}} &= -\grad p\parens{\xx,t} + \mu\grad^2\uu\parens{\xx,t} + \ff\parens{\xx,t}, \label{eq:nes_0}\\
\grad\cdot\uu\parens{\xx,t} &= 0, \label{eq:nes_1}\\
\ff\parens{\xx,t} &= \int_{\SolidRef}\FF\parens{\XX,t}\,\delta\parens{\xx - \bchi\parens{\XX,t}}\,d\XX, \label{eq:lagForce}\\
\frac{\partial \bchi}{\partial t}\parens{\XX,t} &= \int_\Omega \uu\parens{\xx,t}\,\delta\parens{\xx-\bchi\parens{\XX,t}}\,d\xx = \uu\parens{\bchi\parens{\XX,t},t}, \label{eq:motion}
\end{align}
\end{subequations}
in which $\uu\parens{\xx,t}$ and $p\parens{\xx,t}$ are the Eulerian velocity and pressure, respectively, $\ff\parens{\xx,t}$ is the Eulerian force density, $\FF\parens{\XX,t}$ is the Lagrangian force density, which is determined in a manner specfied below, and $\delta\parens{\xx}$ is the Dirac delta function. The fluid density $\rho$ and viscosity $\mu$ are assumed to be constant. \Cref{eq:nes_0,eq:nes_1} are the well known Navier-Stokes equations and hold across the entire computational domain $\Omega$. \Cref{eq:lagForce,eq:motion} couple the Lagrangian and Eulerian variables. The integral in \cref{eq:lagForce} is over the reference configuration of the solid subdomain while that in \cref{eq:motion} is over the entire computational domain.

For the current study, the aortic walls are treated as approximately rigid while the valve leaflets are elastic and deformable. For rigid structures, we use a penalty formulation that is intended to tether the structure in place using the force
\begin{equation}
\FF\parens{\XX,t} = \kappa\parens{\XX - \bchi\parens{\XX,t}},
\end{equation}
in which $\kappa$ is the stiffness parameter \cite{Lee2021}. In practice, we choose $\kappa$ to be the largest stable value permitted by the numerical scheme so that the structure's motion is minimized.

For elastic structures which in this study are the valve leaflets, the response is given by the first Piola-Kirchhoff stress $\PP$, which is determined from a strain-energy function $\Psi\parens{\FFbb}$ via $\PP = \frac{\partial \Psi}{\partial \FFbb},$
in which $\FFbb = \frac{\partial \bchi}{\partial \XX}$ is the deformation gradient tensor. Following the immersed finte element/difference approach of  Vadala-Roth \etal\ \cite{Vadala-Roth2020}, we split the strain energy functional into deviatoric and dilational parts,
$\Psi\parens{\FFbb} = W\parens{\bar{\FFbb}} + U\parens{J},$
in which $J = \text{det}~\FFbb$ is the Jacobian of the deformation tensor and $\bar{\FFbb} = J^{-1/3}\FFbb$. In what follows, we choose the dilational part of the energy to be
\begin{equation}
U\parens{J} = \frac{\kappa_{\text{stab}}}{2}\parens{\log J}^2,
\end{equation}
in which $\kappa_\text{stab}$ is the numerical bulk modulus. The Lagrangian force density is then computed by requiring
\begin{equation}
\int_{\SolidRef} \FF\parens{\XX,t} \cdot\VV\parens{\XX}\,d\XX = -\int_{\SolidRef} \PP\parens{\XX,t}:\grad_\XX \VV\parens{\XX}\,d\XX,
\end{equation}
for all smooth test functions $\VV\parens{\XX}$. 

The leaflets are modeled as a hyperelastic material that follows an exponential neo-Hookean model \cite{Gasser2006,Murdock2018}. The deviatoric strain energy functional for this model is given by
\begin{equation}
W\parens{\bar{\FFbb}} = C_{10}\parens{e^{C_{01}\parens{\bar{I}_1 - 3}} - 1},
\end{equation}
in which $C_{10}$ and $C_{01}$ are material constants and $\bar{I}_1$ is the first deviatoric strain invariant of the modified right Cauchy-Green tensor, $\bar{I}_1 = \text{tr}\parens{\bar{\mathbb{C}}}$, which is defined in terms of the modified deformation tensor, $\bar{\mathbb{C}} = \bar{\FFbb}^\text{T}\bar{\FFbb}$. The material parameter $C_{10}$ is set to be a function of the bound platelet concentration, as described in \cref{sec:mass_deposition}.

\subsection{Boundary Conditions}\label{sec:bdry_conds}
We use reduced order models to determine pressure-flow relationships in the ascending aorta and left ventricular outflow tract (LVOT); see \cref{fig:domains}. These reduced order models are connected to the FSI model through boundary conditions imposed on the fluid \cite{Lee2021a,Griffith2009b}. We use a coupling scheme in which the net flow rate through each of the boundary surfaces serves as an input to the corresponding reduced order model. In turn, each reduced order model determines a pressure that is prescribed on the corresponding boundary surface. The net flow rate through the LVOT boundary is the integral of the vertical component of the velocity over the portion of the bottom boundary of the computation domain between the vesesl walls. The net flow rate for the aortic boundary is defined in a similar way. On the remainder of the computational domain's boundary, we use zero velocity boundary conditions.

We use a three-element Windkessel model \cite{Stergiopulos1999} for the downstream reduced order model that models the aortic outflow,
\begin{align}
C\frac{d\PWK}{d t} &= \QAO - \frac{\PWK}{\RP}, \\
\PAO &= \PWK + \QAO\RC,
\end{align}
in which $C$ is the compliance, $\RC$ is the characteristic resistance, $\RP$ is the peripheral resistance, $\PWK$ is the Windkessel pressure, $\QAO$ is the computed volumetric flow rate at the outlet of the ascending aorta, and $\PAO$ is pressure at the outlet of the ascending aorta which is then prescribed as a boundary condition for the fluid.

For the upstream model to model the inflow from the heart, we employ a time-dependent elastance-based left heart model \cite{Stergiopulos1999},
\begin{align}
\frac{d\parens{\CLA\PLA}}{dt} &= \QVEIN - \QMV, \\
\frac{d\parens{\CLV\PLV}}{dt} &= \QMV - \QLVOT, \\
\PLVOT &= \PLV - \QLVOT\RLVOT, \\
\QMV &= \left\lbrace\begin{array}{ll}
0, &\mbox{ if }\quad \PLA \leq \PLV, \\
\frac{\PLA - \PLV}{\RMV}, &\mbox{ if }\quad \PLA > \PLV,
\end{array}\right.
\end{align}
in which $\CLA$ and $\CLV$ are the time-dependent compliances of the left atrium and left ventricle, respectively. $\RLVOT$ and $\RMV$ are the resistances of the LVOT and mitral valve, the latter of which is modeled as a diode. $\PLA$, $\PLV$, and $\PLVOT$ are the left atrial, left ventricular, and the LVOT pressures, and $\QVEIN$, $\QMV$, and $\QLVOT$ are the volumetric flow rates of the pulmonary veins, mitral valve, and LVOT. In this model, $\QVEIN$ is prescribed as a constant inflow rate into the left atrium. $\QLVOT$ is the computed flow rate at the inlet of the computational domain. $\PLVOT$ is then prescribed as a boundary condition for the momentum \cref{eq:nes_0}. We determine the time-dependent compliances $C(t)$ from specified elastance functions $E(t)$ via $C(t) = 1 / E(t)$. We use the ``two-Hill" elastance waveform given by Mynard \etal\cite{Mynard2012},
\begin{align}
E\parens{t} = \parens{E_\text{max} - E_\text{min}}\alpha\parens{t} + E_\text{min}, \\
\alpha\parens{t} = \frac{k \frac{g_1}{1 + g_1}\frac{1}{1 + g_2}}{\max\parens{\frac{g_1}{1+g_1},\frac{1}{1+g_2}}}, \\
g_i = \left(\frac{t}{\tau_i}\right)^{m_i}.
\end{align}
We use the elastance parameters in $E\parens{t}$ for the left atrium from Mynard \etal\ \cite{Mynard2012}. The remaining parameters are fit to experimental measurements of human aortic pressures $\PAO$ and aortic flow rates $\QAO$ from Murgo \etal\ \cite{Murgo1980} by taking the experimental measurements of $\QAO$ as input to the Windkessel model, and comparing the resulting values of $\PAO$ to its experimental values. We calculate the best-fit parameters to data from Murgo et al.\ \cite{Murgo1980} for a ``Type A" beat for the upstream model. The downstream model is fit using the corresponding downstream data from Murgo et al.\ \cite{Murgo1980}. The fits were created using MATLAB's \verb|fmincon|, a nonlinear optimization tool.

\subsection{Mass Deposition Model}\label{sec:mass_deposition}
We couple the FSI model to a mass deposition model that includes a fluid-phase $\cf\parens{\xx,t}$ concentration measured per unit volume and a surface-bound $\Cb\parens{\XX,t}$ concentration field measured per unit reference area. Although this model does not include the cellular and biochemical interactions describing thrombosis, it does include fields which we view as platelet populations, and the conversion of fluid-phase platelets in $\cf$ to surface-bound platelets in $\Cb$ as platelet adhesion. The fluid-phase species diffuses and advects with the local fluid velocity in the interior fluid domain $\Interior$ and can be converted into the surface-bound species along the boundary $\Gamma_t\subset\InteriorBdry = \partial\Interior$. In the results in \cref{sec:RESULTS}, $\Gamma_t$ is the downstream side of one or both of the leaflets. The surface-bound species moves with the structure and can dissociate to become the fluid-phase species. The model equations are
\begin{subequations}\label{eq:bio_model}
\begin{alignat}{2}
\frac{\partial \cf\parens{\xx,t}}{\partial t} + \uu\parens{\xx,t}\cdot\grad \cf\parens{\xx,t} =& D\grad^2 \cf\parens{\xx,t} ,&& \xx\in\Interior, \\
\frac{\partial \cf\parens{\xx,t}}{\partial \nn} =& 0 ,&& \xx\in\InteriorBdry\setminus\Gamma_t, \label{eq:bdry_neumann}\\
-D\frac{\partial \cf\parens{\xx,t}}{\partial \nn} =& \kon\parens{\Cbmax - \Cb\parens{\bchi\parens{\XX,t},t}}J_\text{s}\cf\parens{\xx,t} &&\nonumber\\
& - \koff\Cb\parens{\bchi\parens{\XX,t},t}J_\text{s} ,&& \xx\in\Gamma_t, \label{eq:bdry_robin}\\
\frac{\partial \Cb\parens{\XX,t}}{\partial t} =& \kon\parens{\Cbmax - \Cb\parens{\XX,t}}\cf\parens{\bchi\parens{\XX,t},t} &&\nonumber\\
&- \koff\Cb\parens{\XX,t} ,&& \XX\in\Gamma_0, \label{eq:bdry_ode}
\end{alignat}
\end{subequations}
in which $D$ is the diffusion coefficient, $\kon$ and $\koff$ are the reaction rates for adhesion and dissociation, respectively, $\Cbmax$ is the carrying capacity of $\Cb$ per unit undeformed area along the boundary $\Gamma_0$, and $J_\text{s} = \frac{dA}{da}$ is the surface Jacobian, which is the ratio of reference and current areas. The first term on the right hand side of \cref{eq:bdry_robin} gives the rate of binding of fluid-phase platelets with concentration $\cf$ to the valve leaflet where the available binding sites have surface density $\parens{\Cbmax - \Cb\parens{\bchi\parens{\XX,t},t}}\Js$ with respect to the current leaflet configuration. The second term gives the rate at which absorbed platelets with surface density $\Cb\parens{\bchi\parens{\XX,t},t}\Js$ detach from the leaflet.

To model the effect of thrombosis over the valve leaflets, we set the stiffness coefficient of the leaflets $C_{10}$ to be a function of the surface concentration $\Cb\parens{\XX,t}$. Because $\Cb\parens{\XX,t}$ is defined only on the surface of the leaflet, we use a harmonic interpolation procedure to extend the surface concentration into the interior of the leaflet, where the Lagrangian forces are calculated. Specifically, we solve Laplace's equation
\begin{align}
\grad^2 \Cb^{\text{in}}\parens{\XX,t} &= 0, \quad \XX\in \Omega^{\text{leaf}}_0, \\
\Cb^{\text{in}}\parens{\XX,t} &= \left\lbrace \begin{array}{ll}
\Cb\parens{\XX,t}, &\mbox{ if } \XX\in\Gamma_0, \\
0, &\mbox{ otherwise},
\end{array}\right.
\end{align}
in which $\Omega^{\text{leaf}}_0$ is the leaflet domain in the initial configuration. Having found $\Cb^\text{in}\parens{\XX,t}$, we then set the stiffness of the leaflet to be
\begin{equation}\label{eq:stiffness}
C_{10}\parens{\XX,t} = C_{10}^\text{base}\parens{\frac{\beta + 1}{2} - \frac{\beta - 1}{2}\cos\parens{\frac{\pi \Cb^\text{in}\parens{\XX,t}}{\Cbmax}}},
\end{equation}
in which $C_{10}^\text{base}$ is the stiffness with no accumulation and $\beta C_{10}^\text{base}$ is the maximum stiffness.

The parameters of the mass deposition model are chosen so that the reactions occur on a similar time scale as the fluid-structure interactions. These values are several orders of magnitude larger than those used in a similar thrombosis model as described previously \cite{Leiderman2011, Fogelson2008}. Use of physiologically relevant reaction rates would require performing simulations over thousands of computational cycles, which is currently not feasible. We are actively working on a temporal multiscale method to meet this challenge and allow use of realistic reaction rates.

\section{Computational Models and Numerical Methods}\label{sec:DISCRETIZATION}
The model is implemented in IBAMR, which provides implementations of the immersed boundary method and several of its extensions along with support for adaptive mesh refinement \cite{ibamr}. IBAMR utilizes libMesh for the finite element representation of the structural deformations \cite{Kirk2006} and PETSc for linear solvers \cite{petsc_user_ref,petsc_web_page,petsc_efficient}. Support for structured adaptive mesh refinement is provided by SAMRAI \cite{Hornung2002}. While the model can be naturally extended to three spatial dimensions, we describe the numerical implementation and results in two spatial dimensions.
\subsection{Imaged Model and Mesh Generation} \label{sec:mesh}
Our two-dimensional aortic root geometry is informed by a three-dimensional patient-specific aortic root model based on pre-procedural computed tomography (CT) image data of a female patient preparing for aortic valve replacement at UNC Medical Center. The images used in this study were obtained under a protocol approved by the UNC Institutional Review Board (study number 18-0202). The CT scan was performed using a Siemens SOMATOM Definition CT Scanner with an image resolution of $512\times 512\times 226$ and a voxel size of $\SI{0.441}{\milli\metre}\times \SI{0.441}{\milli\metre}\times \SI{0.6}{\milli\metre}$. The patient images are segmented by a semi-automated method in ITK-SNAP \cite{itk-snap}, which implements an active contour model that minimizes an energy functional of voxel intensities \cite{Yushkevich2006}. The aortic root measures $\SI{28}{\milli\metre}$ in diameter and $\SI{7.68}{\centi\metre}$ in length, and the thickness of the aortic wall is $\SI{1.0}{\milli\metre}$. The inflow boundary of the model is truncated at the LVOT, and the outflow boundary of the model is truncated downstream of the aortic valve before the first arterial bifurcation. Artificial circular extensions are added at both boundaries using SOLIDWORKS (Dassault Syst\`{e}mes SOLIDWORKS Corporation, Waltham, MA, USA) to simplify the application of boundary conditions to the computational model. The radius of the vessel at both truncations is $\SI{21}{\milli\metre}$. Idealized aortic valve replacement leaflets with a thickness of $\SI{0.7}{\milli\metre}$ are created based on the measurements from Sahasakul \etal\  \cite{Sahasakul1988} and trimmed to fit within the reconstructed aortic root in SOLIDWORKS. To derive our two-dimensional aortic root geometry from the three-dimensional model, we extract a slice through the diameter of the aorta using Coreform Cubit (Computational Simulation Software, LLC, American Fork, UT, USA), which is a software application based on CUBIT from Sandia National Laboratory. We then use Cubit to smooth the angles in both the aortic root and leaflet surfaces and to generate structural meshes consisting of triangular elements.
\subsection{Fluid-Structure Interaction}\label{sec:fsi}
The fluid equations \eqref{eq:nes_0} and \eqref{eq:nes_1} are solved using a second-order Cartesian staggered-grid finite difference method. The nonlinear term is approximated using a piecewise parabolic method \cite{Rider2007}. The resulting saddle point system is solved using GMRES with a projection method as a preconditioner \cite{Griffith2009a}.

The solid subdomain $\Solid$ is discretized using $\mathcal{C}^0$ finite elements. A triangulation $\mathcal{T}_h$ of the structure is constructed. The size of each element in the triangulation is chosen so that there is approximately one node per Cartesian grid cell. On $\mathcal{T}_h$, we define Lagrangian basis functions $\left\{ \phi_l\parens{\XX}\right\}_{l=1}^m$, in which $m$ is the total number of nodes in the triangulation. We approximate the structural deformation and force using the basis functions via
\begin{align}
\bchi\parens{\XX,t} = \sum_{l = 1}^m\bchi_l\parens{t}\,\phi_l\parens{\XX}, \\
\FF\parens{\XX,t} = \sum_{l = 1}^m\FF_l\parens{t}\,\phi_l\parens{\XX}.
\end{align}
Coupling between the fluid and structure is mediated using regularized delta functions in \cref{eq:lagForce,eq:motion}. Recently, Lee and Griffith \cite{Lee2021} suggested using delta functions with smaller support for structures in shear driven regimes. Therefore, in this work, we use the three-point $B$-spline kernel for the flexible valve leaflets, and a two-point piecewise linear kernel for the nearly rigid walls of the aortic root.

\subsection{Mass Deposition Model}
The fluid phase concentration field is approximated using a hybrid semi-Lagrangian cut-cell method \cite{Barrett2022}. For brevity, we omit the details and only highlight the changes of the discretization. To summarize, we split \cref{eq:bio_model} into an advection step,
\begin{equation}\label{eq:advection}
\frac{\partial \cf\parens{\xx,t}}{\partial t} + \uu\parens{\xx,t}\cdot\grad\cf\parens{\xx,t} = 0,\quad \xx\in\Interior,
\end{equation}
and a diffusion step,
\begin{equation}\label{eq:diffusion}
\frac{\partial \cf\parens{\xx,t}}{\partial t} = D\grad^2\cf\parens{\xx,t}, \quad\xx\in\Omega_t^{f^-},
\end{equation}
along with the boundary conditions \cref{eq:bdry_neumann,eq:bdry_robin} and the surface concentration \cref{eq:bdry_ode}. During the diffusion step, the domain $\Interior$ is assumed to be fixed. The advective step is treated with a semi-Lagrangian method using polyharmonic splines to reconstruct the function $\cf$. The diffusion step is treated with a cut-cell finite volume method. The surface concentration $\Cb\parens{\XX,t}$ is solved for by extrapolating the fluid-phase field $\cf\parens{\xx,t}$ to the boundary and approximating the ODE in \cref{eq:bdry_ode}.

\subsubsection{Diffusion}\label{sec:diffusion}
To approximate the diffusion step in \cref{eq:diffusion}, we employ a cut-cell finite volume method in which the domain $\Omega_t^{f^-}$ is considered fixed for the duration of this step. Integrating \cref{eq:diffusion} over a grid cell $\cc_{i,j}$ that is entirely or partially interior to $\Omega_t^{f^-}$ and dividing by that cell's volume, we obtain
\begin{equation}\label{eq:diffusion_int}
\frac{1}{\size{\cc_{i,j}\cap\Omega_t^{f^-}}}\int_{\cc_{i,j}\cap\Omega_t^{f^-}} \frac{\partial \cf\parens{\xx,t}}{\partial t} \text{d}\xx = \frac{1}{\size{\cc_{i,j}\cap\Omega_t^{f^-}}}\int_{\cc_{i,j}\cap\Omega_t^{f^-}} D\Lap\cf\parens{\xx,t} \text{d}\xx.
\end{equation}
We define $C_{\text{f},i,j}$ as the cell average of $\cf\parens{\xx,t}$ in the cell $\cc_{i,j}\cap\Omega_t^{f^-}$. Replacing the cell average in the left hand side of \cref{eq:diffusion_int} and employing the divergence theorem on the right hand side, we obtain
\begin{equation}\label{eq:diffusion_div}
\frac{\text{d}C_{\text{f},i,j}}{\text{d} t} = \frac{1}{\size{\cc_{i,j}\cap\Omega_t^{f^-}}}\int_{\partial\parens{\cc_{i,j}\cap\Omega_t^{f^-}}} D\frac{\partial \cf\parens{\xx,t}}{\partial \nn}\cdot \text{d}\AA.
\end{equation}
The integral in \cref{eq:diffusion_div} consists of two parts, an integral over the boundary $\InteriorBdry$ that is interior to cell $\cc_{i,j}$ and an integral over the portion of the boundary of the cell $\cc_{i,j}$ that is interior to $\Interior$. The first type consists of an integral over the physical boundary and using the provided boundary conditions in \cref{eq:bdry_neumann,eq:bdry_robin}, can be computed using techniques described in the next section. The second integral is discretized using second order finite differences. This discretization requires the computation of the cut cell volume $\size{\cc_{i,j}\cap\Omega_t^{f^-}}$, which is described in \cref{sec:cut_cell_geom}.

\subsubsection{Surface Reactions}\label{sec:surface}
Along part of the surface $\InteriorBdry$, we allow for binding of the fluid-phase species to the boundary and for unbinding of the surface-bound species into the fluid, as described by \cref{eq:bdry_robin,eq:bdry_ode}. We extract a boundary mesh represented by $C^0$ elements from the volumetric leaflet mesh as described in \cref{sec:fsi}. We maintain a representation of both the surface concentration $\Cb\parens{\XX,t}$ per unit reference area and the fluid concentration $\cf\parens{\XX,t}$ per unit volume restricted to the boundary. These values are represented using Lagrangian basis functions $\left\{\psi_l\parens{\XX}\right\}_{l=1}^{n_\text{bd}}$ in which $n_\text{bd}$ is the number of nodes of the boundary mesh. We note these are the same basis functions used for the structural deformation, but restricted to the surface. The concentrations along the boundary are accordingly
\begin{align}
\Cb\parens{\XX,t} = \sum_{l = 1}^{n_\text{bd}} \Cb^l\parens{t}\psi_l\parens{\XX}, \\
\cf\parens{\XX,t} = \sum_{l = 1}^{n_\text{bd}} \cf^l\parens{t}\psi_l\parens{\XX}.
\end{align}
The values $\cf^l$ are found by using a radial basis function interpolant as described in \cref{sec:reconstructions} to extrapolate the values of $\cf\parens{\xx,t}$ to the surface nodes. The nodal values $\Cb^l$ are found by solving the ODE in \cref{eq:bdry_ode} using a two stage Runge Kutta method.

This finite element representation allows for easy evaluations of the flux defined in \cref{eq:bdry_robin} from the boundary to the fluid. To evaluate this flux, we require the value of the Jacobian $J_\text{s} = \frac{d\aa}{d\AA}$ that converts areas in the reference configuration to areas in the current configuration. Because we are using a $C^0$ representation of the surface, the Jacobian is discontinuous at nodes. To obtain a continuous representation, we project $J_\text{s}$ onto the finite element basis \cite{Kolahdouz2020}. In practice, this amounts to computing the Jacobian at quadrature points along the surface.

\subsubsection{Reconstructions}\label{sec:reconstructions}
Both the semi-Lagrangian step and the surface reactions involve reconstructing $\cf\parens{\xx,t}$ at various points $\hat{\xx}$ in the computational domain. The details of the reconstruction procedure depend on where the reconstruction is being performed within this domain. Away from the boundary, we use the four closest grid points to $\hat{\xx}$ to form a bilinear interpolant. If there are too few points to form the bilinear interpolant (e.g., near cut-cells), we use a radial basis function (RBF) interpolant \cite{Flyer2016,Shankar2018}. The RBF interpolant is constructed via a polyharmonic spline
\begin{equation}
q\parens{\xx} = \sum_{j = 1}^k\lambda_j\norm{\xx - \xx_j}^m + \sum_{j = 1}^s \beta_j p_j\parens{\xx},
\end{equation}
in which $m$ is an odd integer and $p_j\parens{\xx}$ form a set of $s$ polynomial basis functions. The total number of points in the stencil $k$ is chosen so that $k = 2m+1$. The points $\xx_j$ are the $k$ closest points to the location $\hat{\xx}$. We find the coefficients $\lambda_j$ and $\beta_j$ by requiring
\begin{subequations}
\begin{alignat}{2}
q\parens{\xx_j} &= f_j &&\mbox{ for } j=1,\ldots,k, \label{eq:interp_conds}\\
\sum_{i= 1}^s \lambda_i p_i\parens{\xx_j} &= 0 &&\mbox{ for } j=1,\ldots,k. \label{eq:orthog_conds}
\end{alignat}
\end{subequations}
\Cref{eq:interp_conds} are the interpolation conditions, and \cref{eq:orthog_conds} are the orthogonality conditions to ensure a unique interpolant. This results in a linear system for the coefficients, which is solved using a QR algorithm. In our computations, we set the integer $m = 3$ and use up to quadratic polynomials.

\subsubsection{Cut Cell Geometries}\label{sec:cut_cell_geom}
In the cut-cell finite volume discretization of \cref{eq:diffusion}, we require the computation of the geometries of cells cut by the boundary $\Gamma^{f^-}$. We denote the node of a Carteisan grid cell by $\xx_{i+\frac{1}{2},j+\frac{1}{2}}$. To find cut cell volumes, we first calculate the signed distance function to the surface at each node $\xx_{i+\frac{1}{2},j+\frac{1}{2}}$. To do this, we first find intersections of the $C^0$ representation of the immersed structure with the background Eulerian grid, and for each element of the immersed structure, we calculate outward facing normals. We note that this requires a consistent traversal (e.g., counter-clockwise) of the structure to ensure a consistent facing normal. Then, for each cell node $\xx_{i+\frac{1}{2},j+\frac{1}{2}}$ of the background grid, we find the projection of the point onto each element and compute its distance from $\xx_{i+\frac{1}{2},j+\frac{1}{2}}$. If multiple minimal distance projections exist, we use the angle weighted average of the projections \cite{Baerentzen2005}. The sign of the distance is computed using the previously computed structure normal. Once we have the signed distances at each cell node, we can compute partial cell volumes. Following Min and Gibou \cite{Min2007}, we compute cell volumes by decomposing the cell into simplices, for which analytic formulas for the volume exist.

\subsection{Time Stepping}
In summary, the steps to advance the solution from time $t^n$ to time $t^{n+1}$ are:
\begin{enumerate}
\item Compute the cut cell geometries.
\item Perform a half step of the diffusion solve, evolving both the fluid-phase and the structure-bound concentration fields.
\item Solve the Navier-Stokes equations and update the position of the immersed structure.
\item Update the cut cell geometries using the new position of the immersed structure.
\item Perform a full step of the semi-Lagrangian method, using the velocities from the Navier-Stokes solve.
\item Perform a half step of the diffusion step, evolving the fluid-phase and surface-bound concentrations.
\end{enumerate}

Our use of an explicit time stepping scheme for several of these steps limits our time step size to resolve the fastest time scale. In this case, the fastest time scale is that of the leaflet elasticity. We determine an empirical scaling relationship between the time step size and the stiffness of the leaflet that maintains numerical stability under increasing leaflet stiffness. Specifically, we choose the time step such that
\begin{equation}
\Delta t = \frac{C_\text{ts}}{\sqrt{\text{max}\parens{C_{10}}}},
\end{equation}
in which $C_\text{ts}$ is chosen to be as large as possible.

\section{Results}\label{sec:RESULTS}
\Cref{tab:parameters} provides the values of all relevant physical and numerical parameters. At $t = 0$, the initial fluid phase concentration $\cf$ is set to be 1 throughout the domain $\Interior$. During the first two cycles, the binding and unbinding coefficient, $\kon$ and $\koff$, are set to zero; afterward they are reset to their non-zero values. We emphasize that the binding and unbinding coefficients and the diffusion coefficient are artificially increased by several orders of magnitude compared to other clotting models \cite{Leiderman2011, Fogelson2008} to ensure that sufficient binding can occur within the duration of the simulation.

\begin{table}
\begin{center}
\caption{Values of the parameters used in the simulation.}\label{tab:parameters}
\begin{tabular}{cc|cc}
\multicolumn{2}{c|}{Structure parameters}         & \multicolumn{2}{c}{Deposition and Fluid parameters} \\
\hline
$\kappa$ & \SI{20.17}{\giga\pascal\per\square\cm} & $\kon$ & \SI{0.03321}{\cubic\cm\per\s\per\platelet} \\
$k_\text{stab}$ & \SI{58.4}{\mega\pascal}         & $\koff$ & \SI{0.01}{\per\s} \\
$C_{01}$ & \SI{3.25}{}                            & $D$ & \SI{0.1}{\square\cm\per\s} \\
$C_{10}^\text{min}$ & \SI{2.264}{\mega\pascal}    & $C_\text{b}^\text{max}$ & \SI{1.41e7}{\platelet\per\square\cm} \\
$\beta$ & varies, between $1-600$                                  & $c_\text{f}^\text{max}$ & \SI{1.5e5}{\platelet\per\cubic\cm} \\
&                                                 & $\rho$ & \SI{1}{\gram\per\cubic\cm} \\
&                                                 & $\mu$ & \SI{0.035}{\gram\per\cm\per\s} \\
\hline
\multicolumn{4}{c}{Boundary model parameters} \\
\hline
$\RP$ & \SI{0.9046}{\mmHg\second\per\milli\liter}   & $\tau_{1,\text{LA}}$ & \SI{0.09789}{\second} \\
$C$ & \SI{1.950}{\milli\liter\per\mmHg}                    & $\tau_{1,\text{LV}}$ & \SI{0.0887}{\second} \\
$\RC$ & \SI{0.042}{\mmHg\second\per\milli\liter}    & $m_{1,\text{LA}}$    & \SI{1.32}{} \\
$\QVEIN$ & \SI{6.2}{\liter\per\minute}              & $m_{1,\text{LV}}$    & \SI{2.404}{} \\
$\QLVOT$ & \SI{0.015}{\mmHg\second\per\milli\liter} & $\tau_{2,\text{LA}}$ & \SI{0.1602}{\second} \\
$\RMV$ & \SI{0.005}{\mmHg\second\per\milli\liter}   & $\tau_{2,\text{LV}}$ & \SI{0.4461}{\second} \\
$E_\text{max,\text{LA}}$       & \SI{0.17}{\mmHg\per\milli\liter}    & $m_{2,\text{LA}}$    & \SI{13.1}{} \\
$E_\text{min,\text{LA}}$       & \SI{0.08}{\mmHg\per\milli\liter}    & $m_{2,\text{LV}}$    & \SI{20.952}{} \\
$E_\text{min,\text{LV}}$       & \SI{0.0265}{\mmHg\per\milli\liter}    & $E_\text{max,\text{LV}}$    & \SI{0.16}{\mmHg\per\milli\liter}
\end{tabular}
\end{center}
\end{table}

\subsection{Convergence Study}

\begin{figure}
\psfragscanon
\psfrag{y}{Velocity Magnitude $\left(\SI{}{\cm\per\second}\right)\quad$}
\psfrag{a}{(b)}
\psfrag{b}{(c)}
\psfrag{c}{(d)}
\psfrag{d}{(e)}
\psfrag{Uv}{}
\psfrag{6v}{0}
\psfrag{5v}{23}
\psfrag{4v}{46}
\psfrag{3v}{69}
\psfrag{2v}{92}
\psfrag{1v}{115}
\psfrag{1}{20}
\psfrag{2}{40}
\psfrag{3}{60}
\psfrag{4}{80}
\psfrag{5}{100}
\psfrag{6}{120}
\begin{tabular}{ccc}
\multirow{2}{*}[6em]{\shortstack{\subfigimg[width=0.3\textwidth,pos=ul,hsep=1.5em,vsep=-0.4em]{(a)}{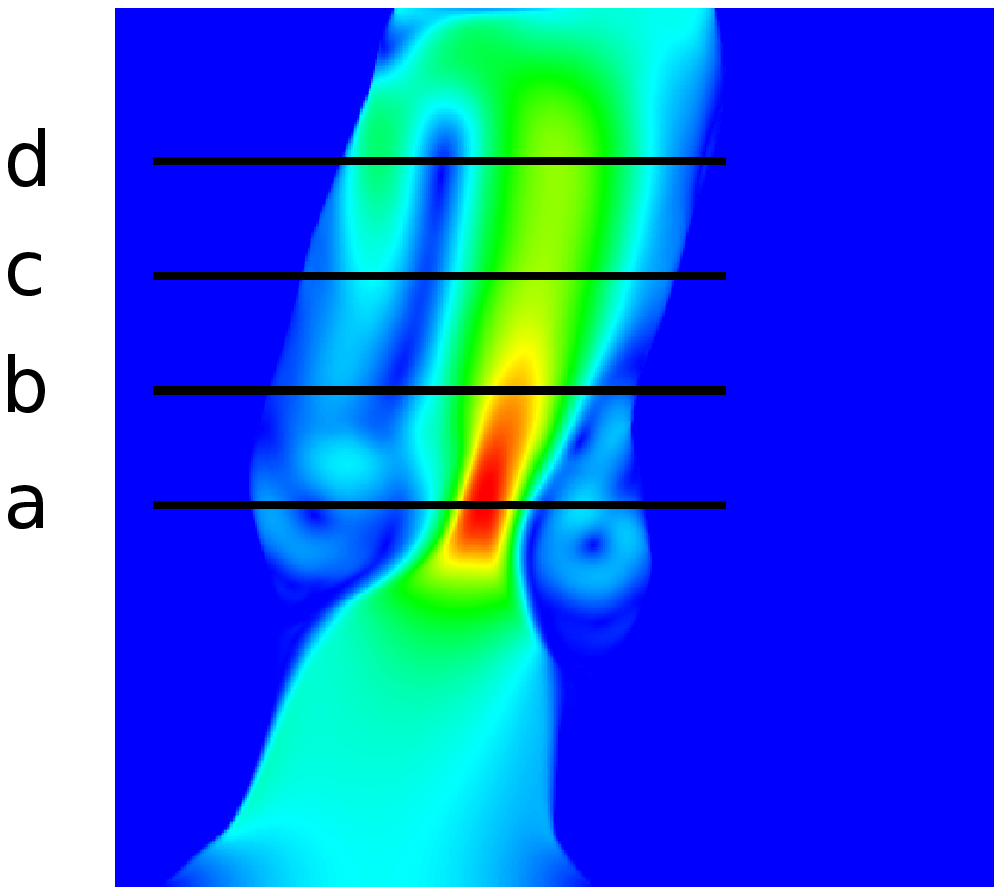} \\ \hspace{0.5em}\includegraphics[width=0.3\textwidth]{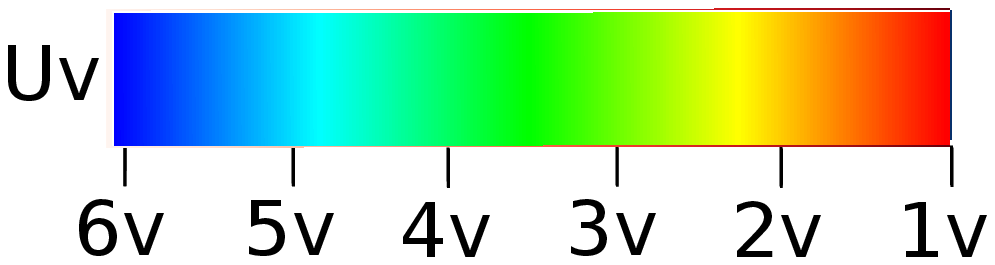}}}
&
\subfigimg[width=0.3\linewidth,pos=ul]{(b)}{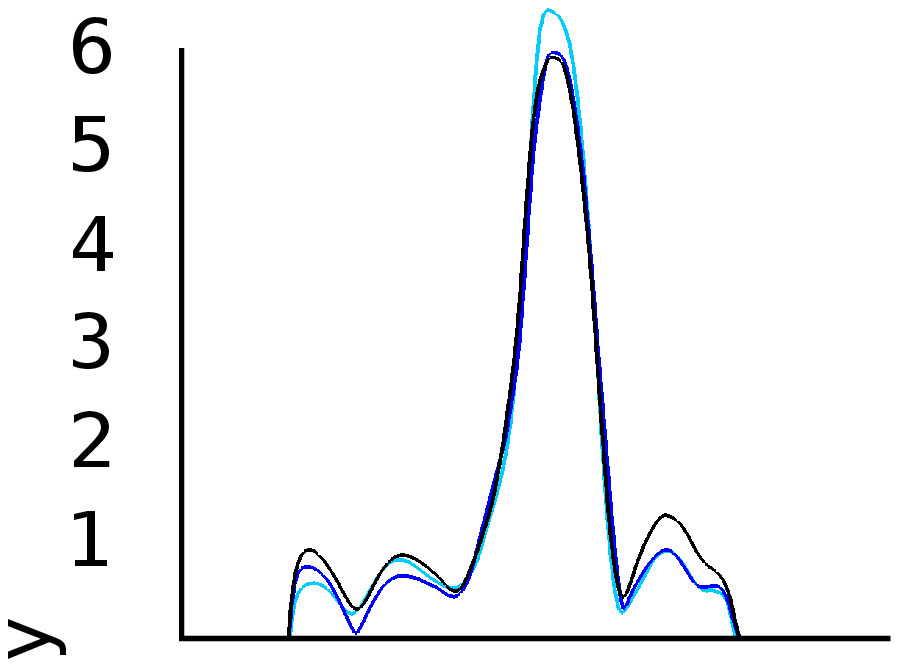} 
&
\subfigimg[width=0.3\linewidth,pos=ul]{(c)}{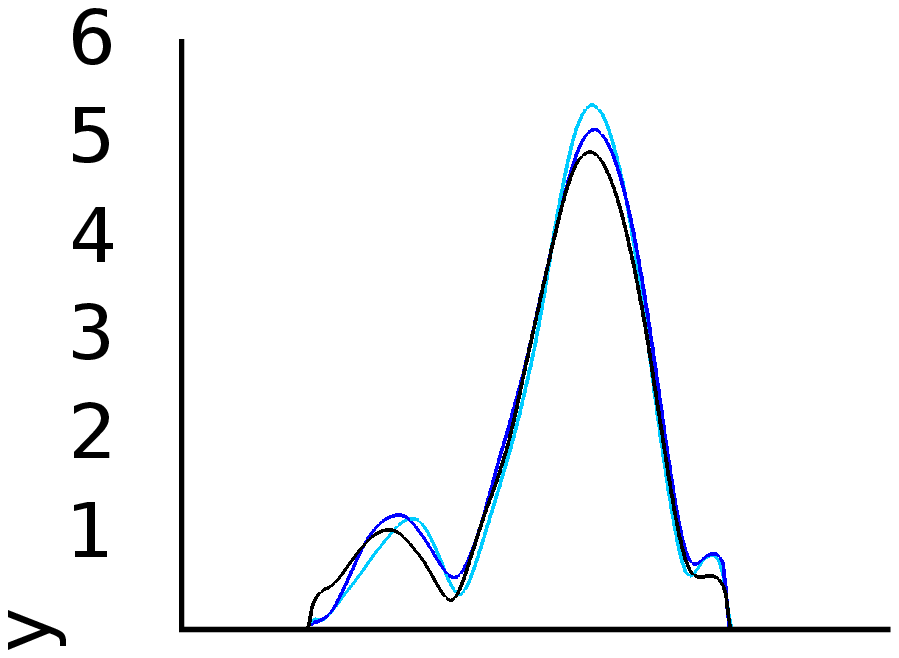} \\
&
\subfigimg[width=0.3\linewidth,pos=ul]{(d)}{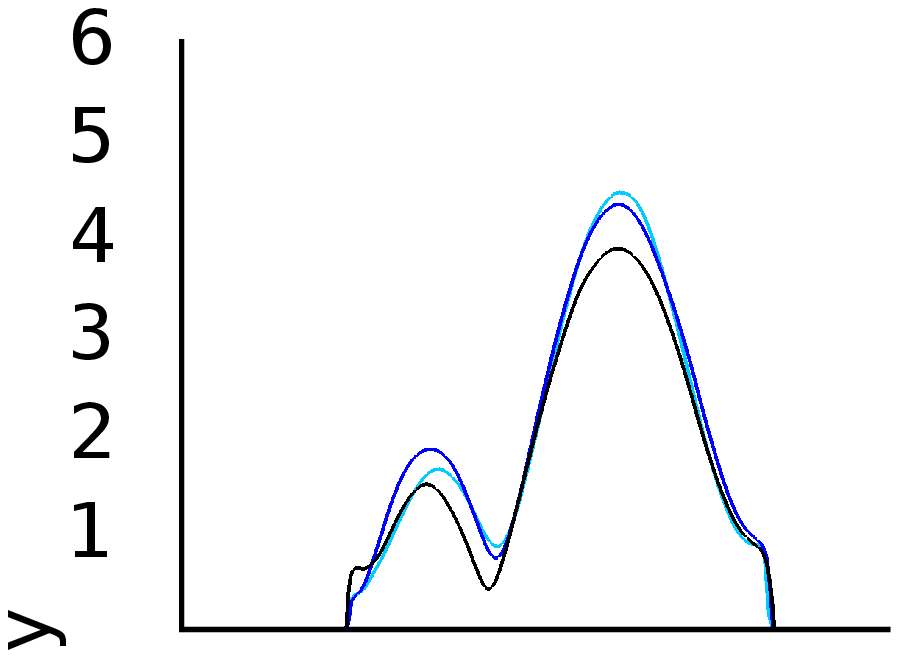} 
&
\subfigimg[width=0.3\linewidth,pos=ul]{(e)}{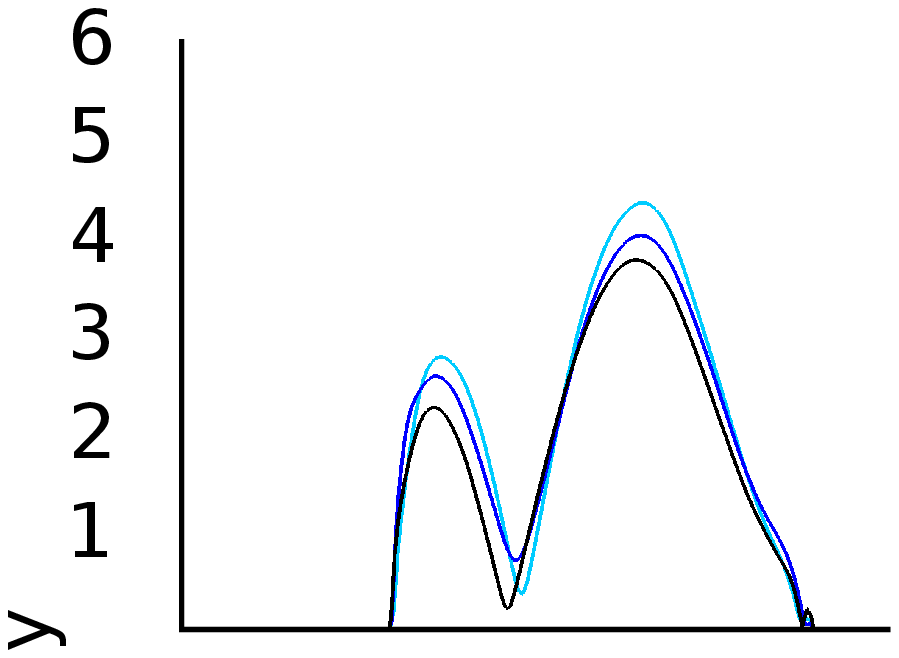}
\end{tabular}
\caption{To assess convergence, we perform simulations at approximately three-quarters systole, and we compute the average velocity from time $0.5$ to time $T = 3.5$. Panel (a) shows the average fluid velocity magnitude across the time interval. Panels (b)-(e) show slices of the average fluid velocity magnitude for three different grid resolutions. The coarsest grid is shown in light blue, the medium grid is shown in blue, and the finest grid is shown in black. We observe consistent values across all grid resolutions tested.}
\label{fig:converge_avg}
\end{figure}

\begin{figure}
    \centering
    \includegraphics[width=0.6\textwidth]{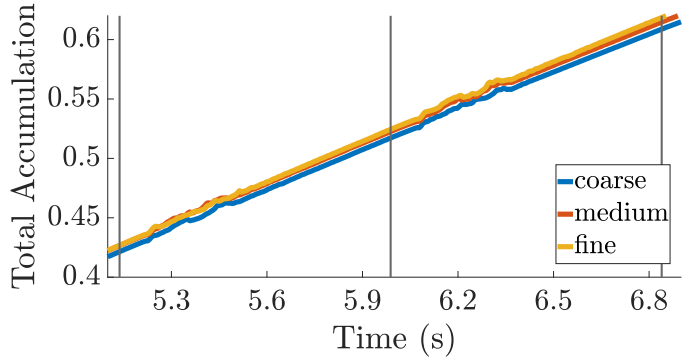}
    \caption{The total accumulation $\int_{\Gamma_0}\Cb\parens{\XX,t}d\XX$ over time for three different mesh sizes. Although we do not expect pointwise convergence of the fluids phase concentration and velocity fields for the turbulent flow regime considered in this study, we do observe grid independence of the bound concentration field.}
    \label{fig:conv_tot_accum}
\end{figure}

The flow regime during peak systole is turbulent, with the largest Reynolds number being approximately 5000. Because of the chaotic nature of the simulation, convergence of the numerical method is not well defined. Small changes in the simulations (e.g. grid spacing and time step size) can lead to large changes in the flow velocities. Further, the fluid-phase and surface concentrations and hence the stiffness of the leaflets is directly affected by the turbulent flow. Therefore, to assess the accuracy of the simulation, we compare the average fluid velocity near peak systole across grid sizes. We modify the model to use a parabolic velocity profile which corresponds to three-quarters systole. We generate three different resolutions of meshes with maximum element edge lengths of 0.32, 0.24, and \SI{0.18}{\mm}, which correspond to 2, 3, and 4 elements across the width of the leaflets, respectively. The background Cartesian grid is refined such that there is approximately one structural mesh node per grid cell. The modified model is then run without accumulation. \Cref{fig:converge_avg} shows the average fluid velocity from time $t = 0.5$ to a final time $T = 3.5$. We observe consistent values across all grid resolutions tested. While convergence is not clear, we expect the average flow velocity to show convergence in the limit as $T\rightarrow\infty$. In the full model, we do observe grid independence of the surface concentration field. \Cref{fig:conv_tot_accum} shows the total bound concentration $\int_{\Gamma_0} \Cb\parens{\XX,t}d\XX$ for all three grid resolutions. For the results presented below, we use the coarse mesh, consisting of two elements across the leaflet.

\subsection{Leaflet Deposition}
\begin{figure}
\begin{center}
\psfragscanon
\psfrag{Us}{\SI{}{\platelet\per\cm}}
\psfrag{1s}{\SI{0.220}{}}
\psfrag{2s}{\SI{0.211}{}}
\psfrag{3s}{\SI{0.201}{}}
\psfrag{4s}{\SI{0.189}{}}
\psfrag{5s}{\SI{0.173}{}}
\psfrag{6s}{\SI{0.150}{}}
\psfrag{Uf}{\SI{}{\platelet\per\cm\squared}}
\psfrag{1f}{\SI{1.000}{}}
\psfrag{2f}{\SI{0.969}{}}
\psfrag{3f}{\SI{0.933}{}}
\psfrag{4f}{\SI{0.888}{}}
\psfrag{5f}{\SI{0.831}{}}
\psfrag{6f}{\SI{0.750}{}}
\psfrag{Uv}{\SI{}{\cm\per\second}}
\psfrag{1v}{\SI{275}{}}
\psfrag{2v}{\SI{202}{}}
\psfrag{3v}{\SI{141}{}}
\psfrag{4v}{\SI{88}{}}
\psfrag{5v}{\SI{41}{}}
\psfrag{6v}{\SI{0.0}{}}
\phantomsubcaption\label{fig:during_sim_right:diastole}
\phantomsubcaption\label{fig:during_sim_right:peak_sys}
\phantomsubcaption\label{fig:during_sim_right:end_sys}
\phantomsubcaption\label{fig:during_sim_right:final}
\begin{tabular}{cccc}
\multicolumn{4}{c}{\includegraphics[width=0.5\linewidth,height=25bp,trim=0 0 0 0]{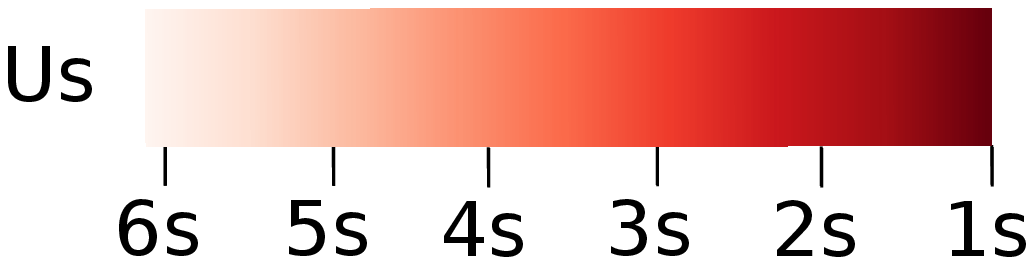}
\includegraphics[width=0.5\linewidth,height=25bp,trim=0 0 0 0]{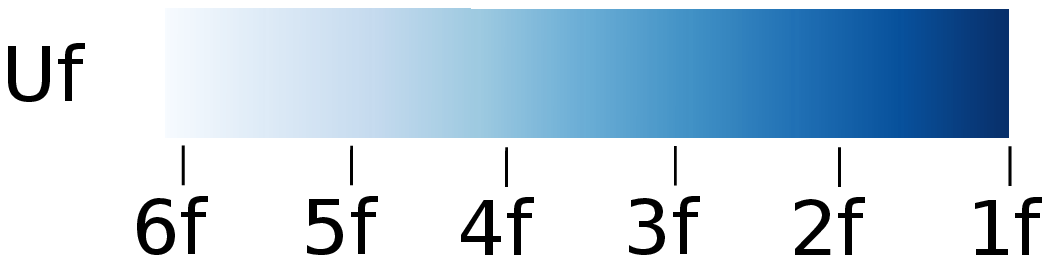}} \\
\subfigimg[width=0.25\linewidth,hsep=0.5em,vsep=1em,pos=ul]{(a)}{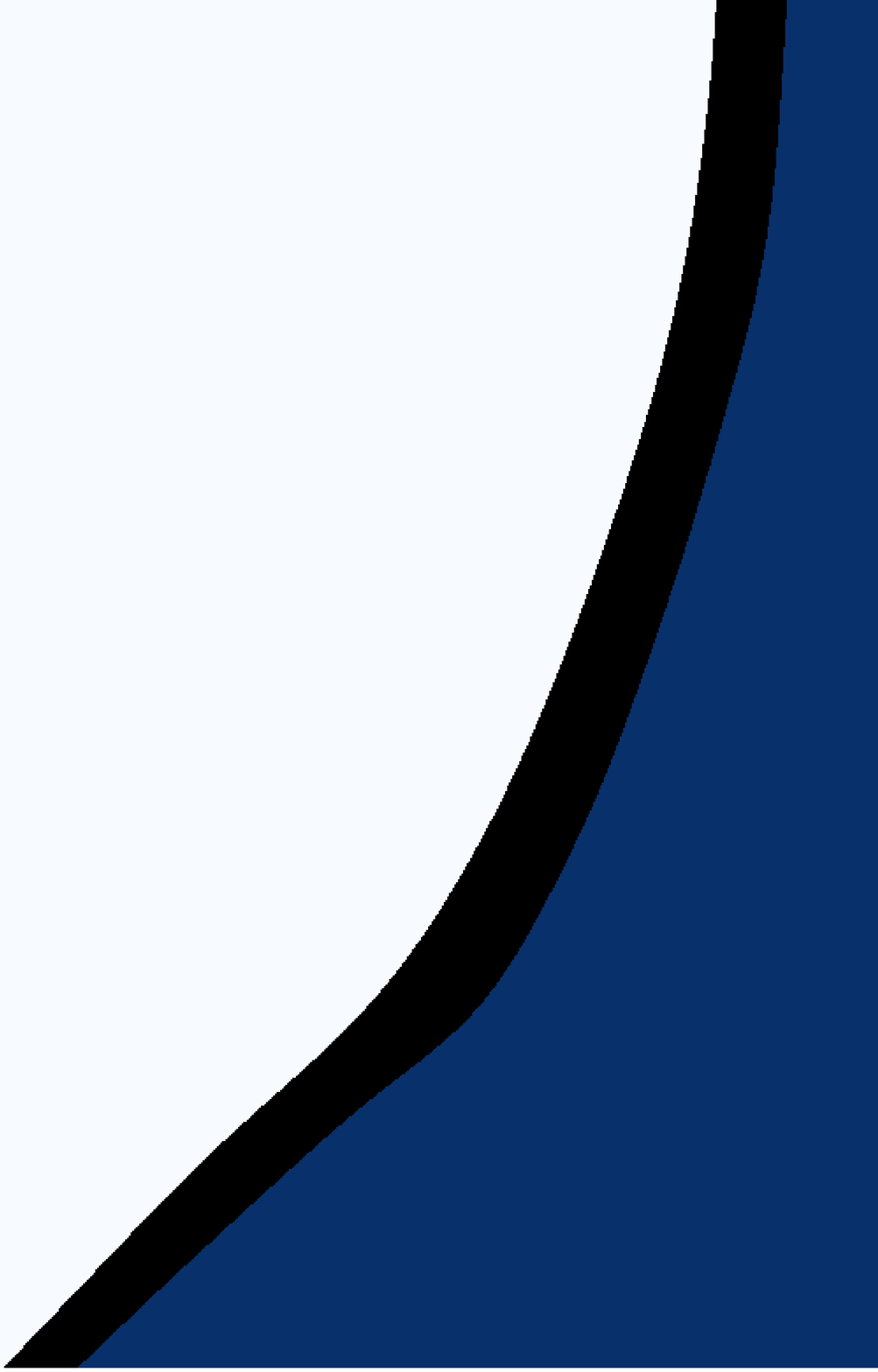}
&
\subfigimg[width=0.25\linewidth,hsep=0.5em,vsep=1em,pos=ul]{(b)}{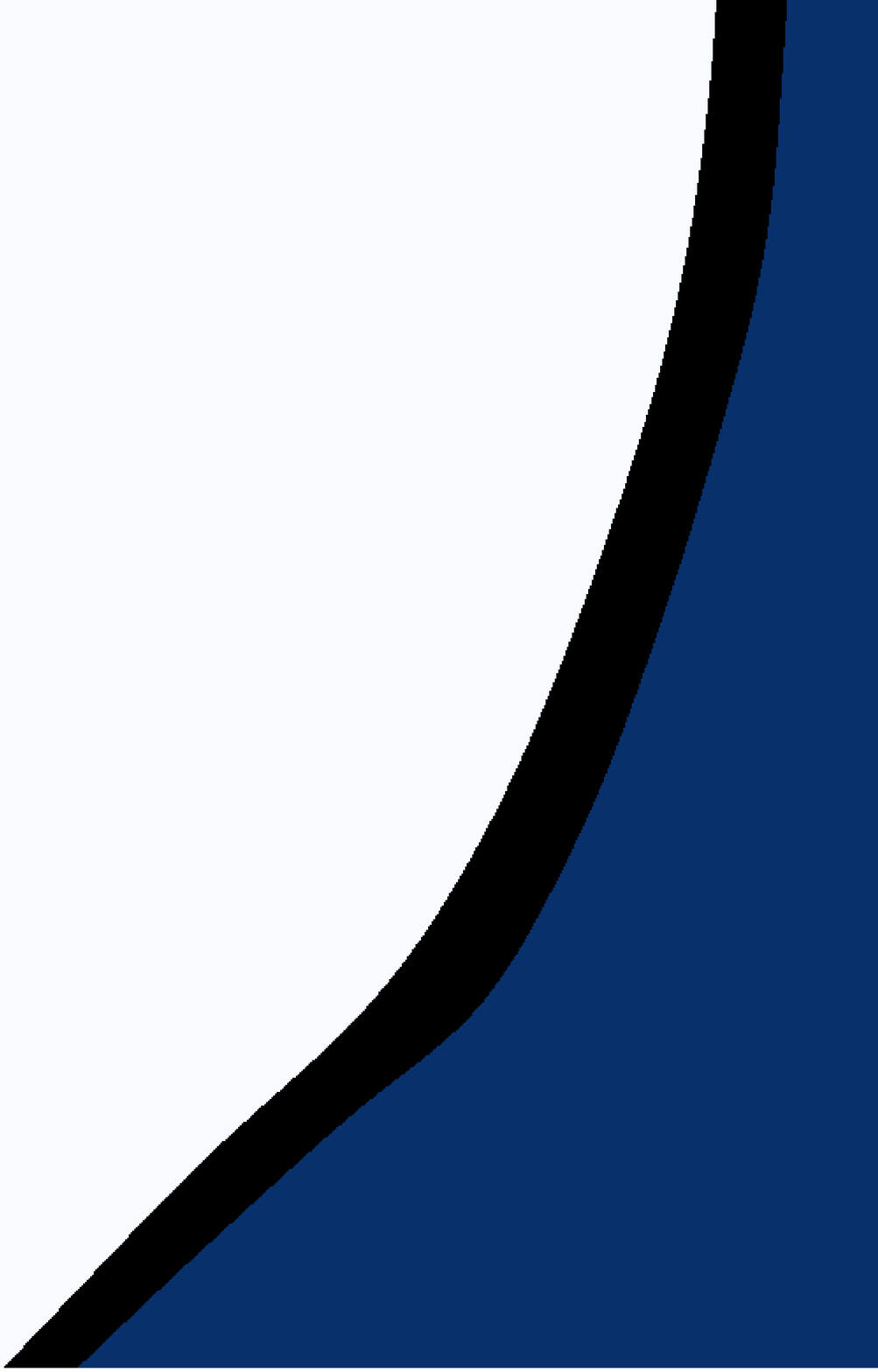}
&
\subfigimg[width=0.25\linewidth,hsep=0.5em,vsep=1em,pos=ul]{(c)}{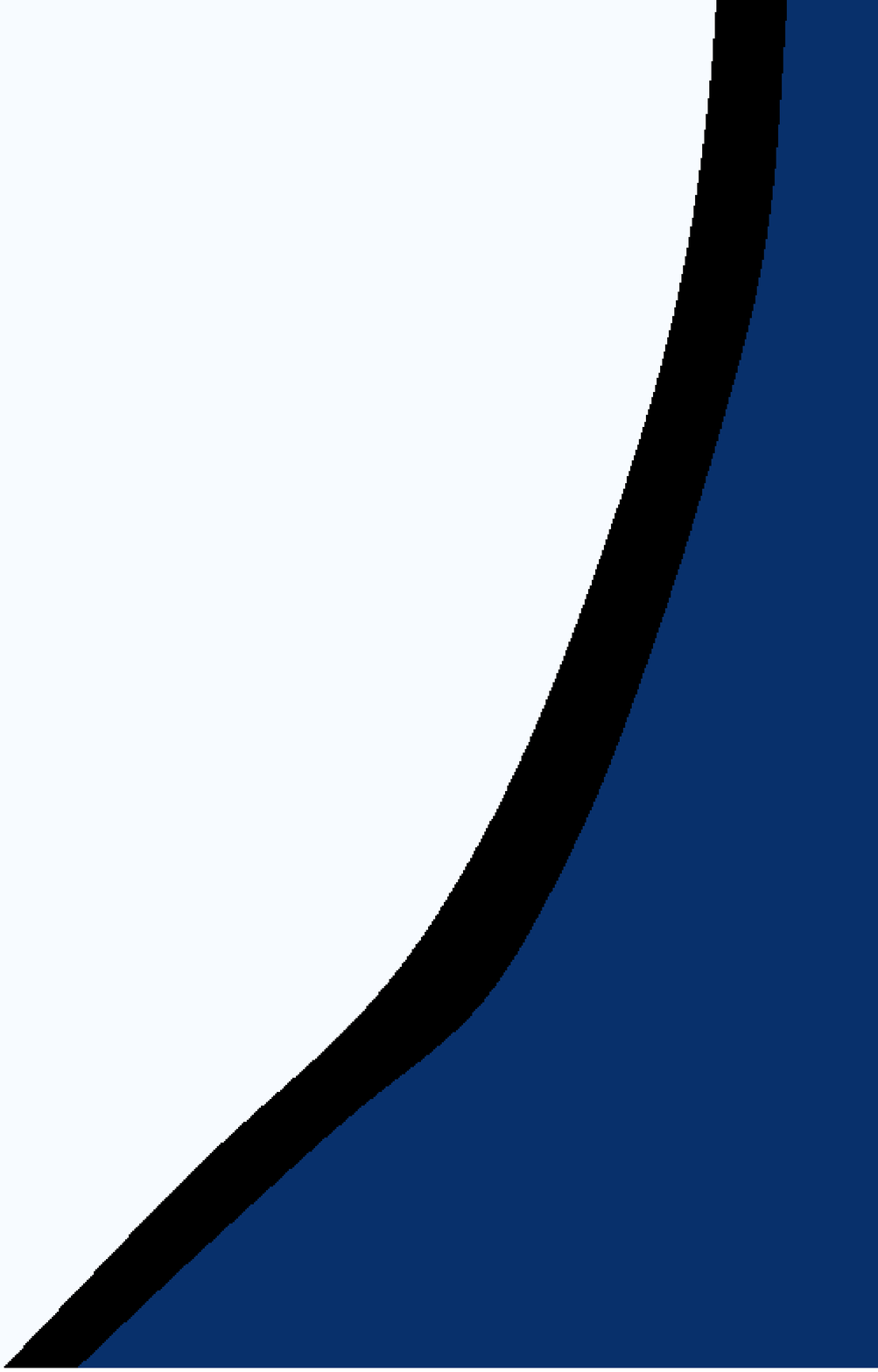}
&
\subfigimg[width=0.25\linewidth,hsep=0.5em,vsep=1em,pos=ul]{(d)}{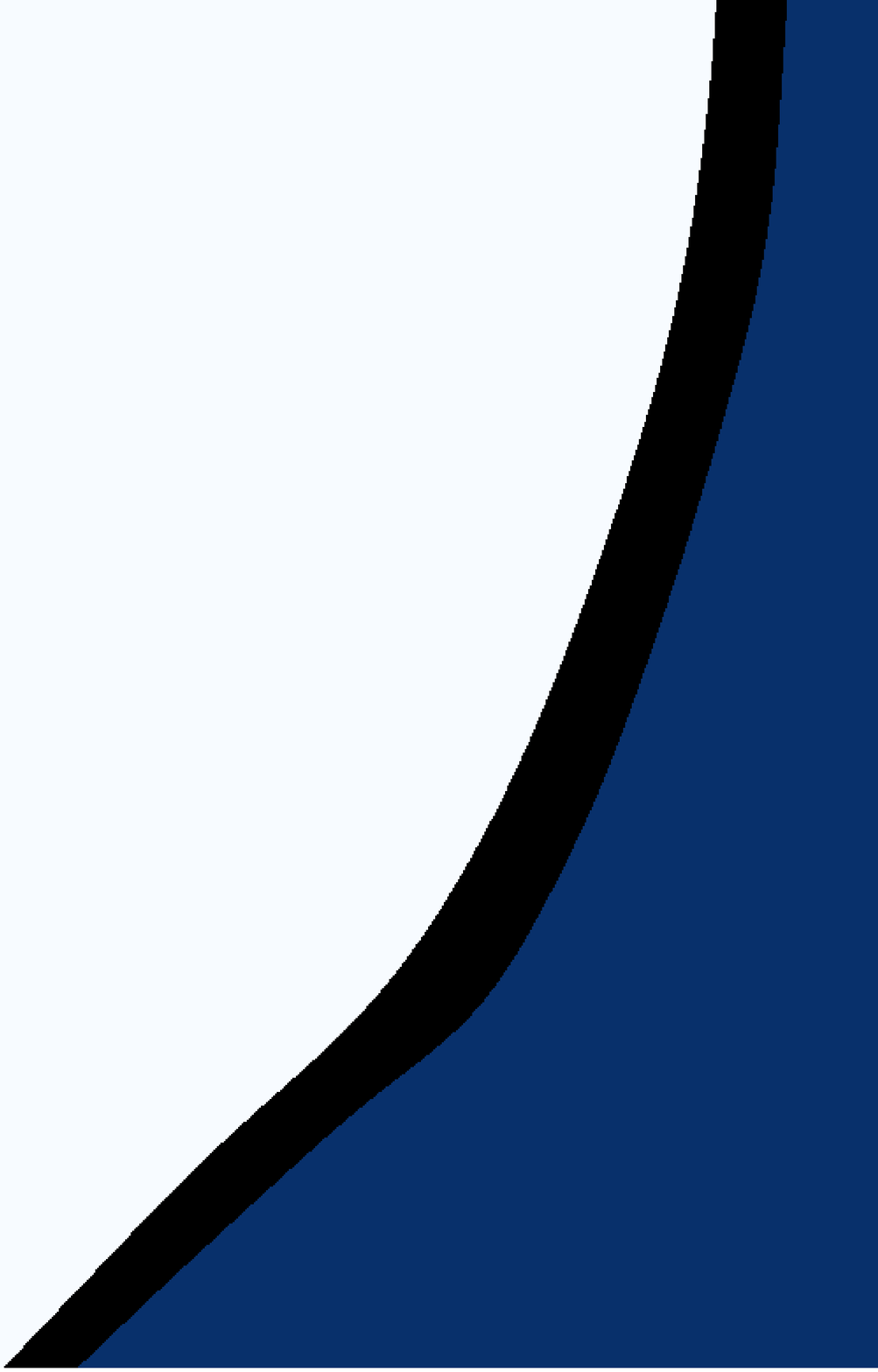} \\
\includegraphics[width=0.25\linewidth]{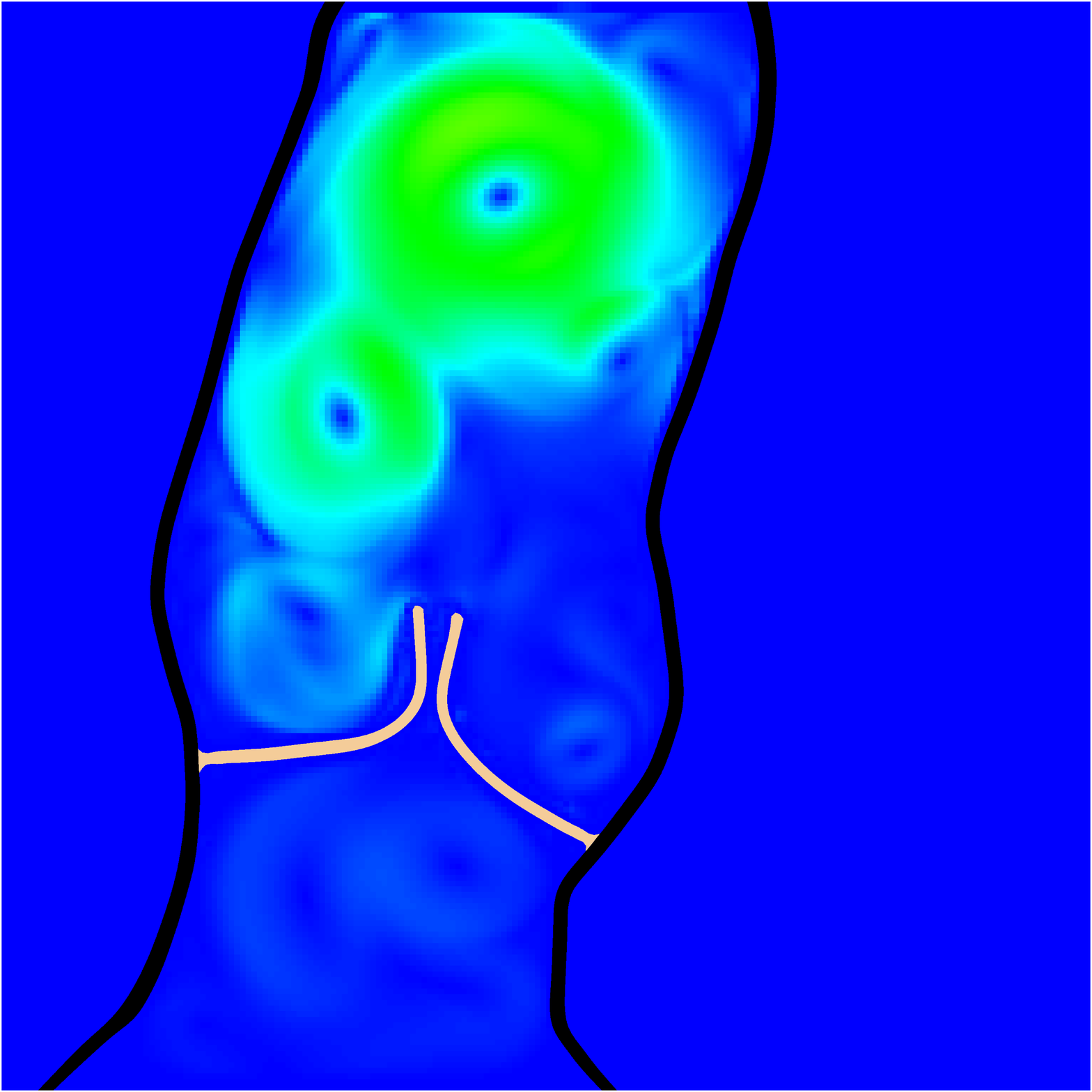}
&
\includegraphics[width=0.25\linewidth]{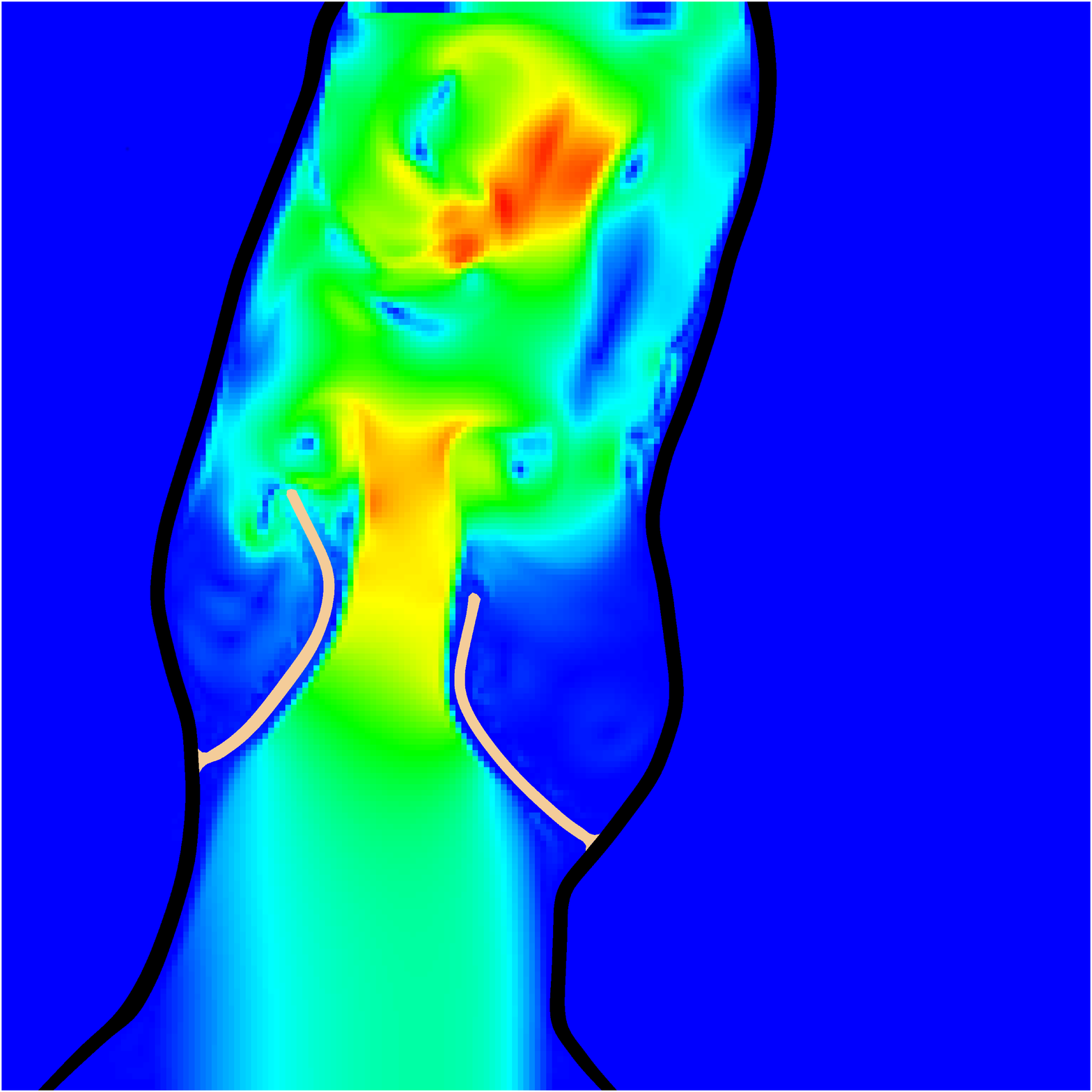}
&
\includegraphics[width=0.25\linewidth]{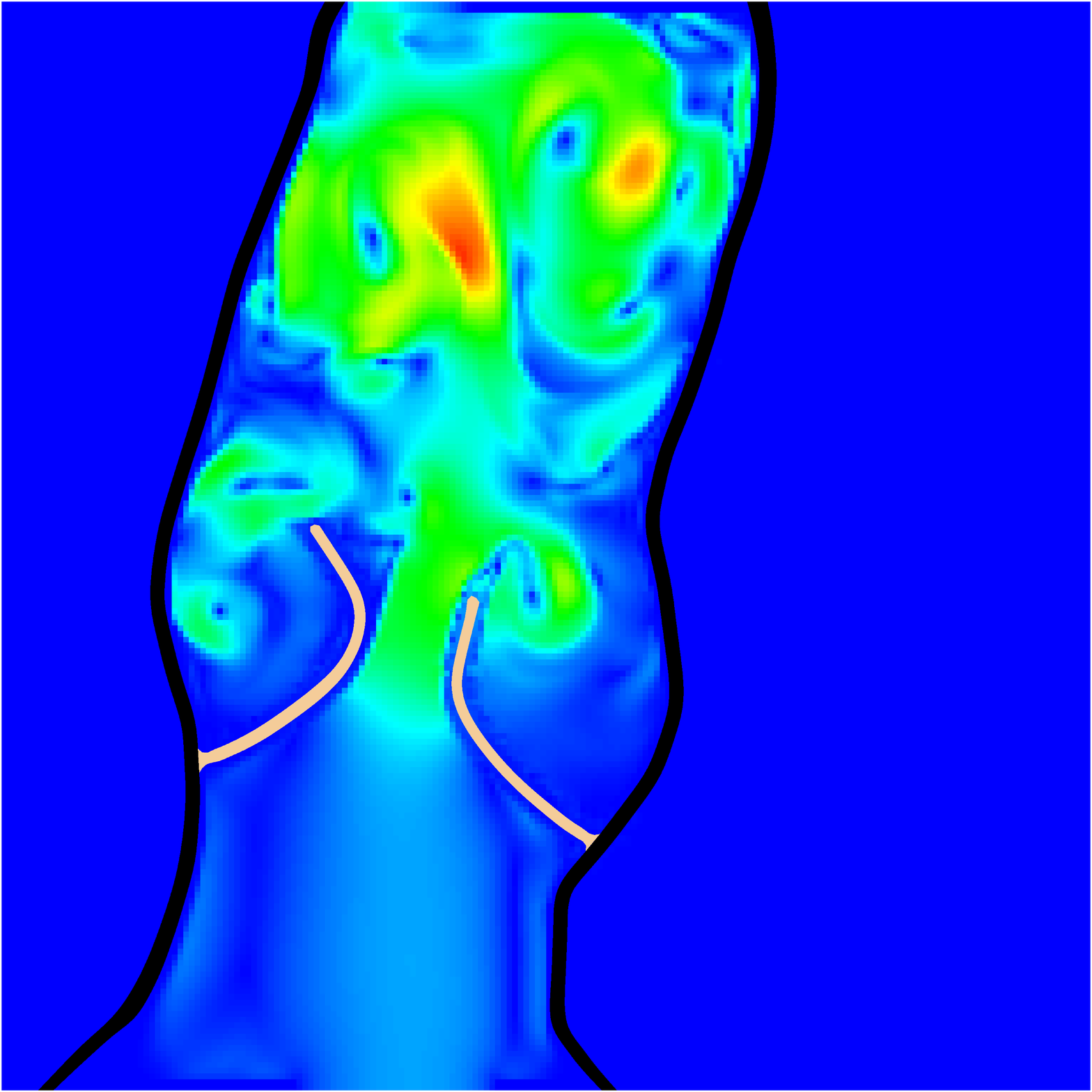}
&
\includegraphics[width=0.25\linewidth]{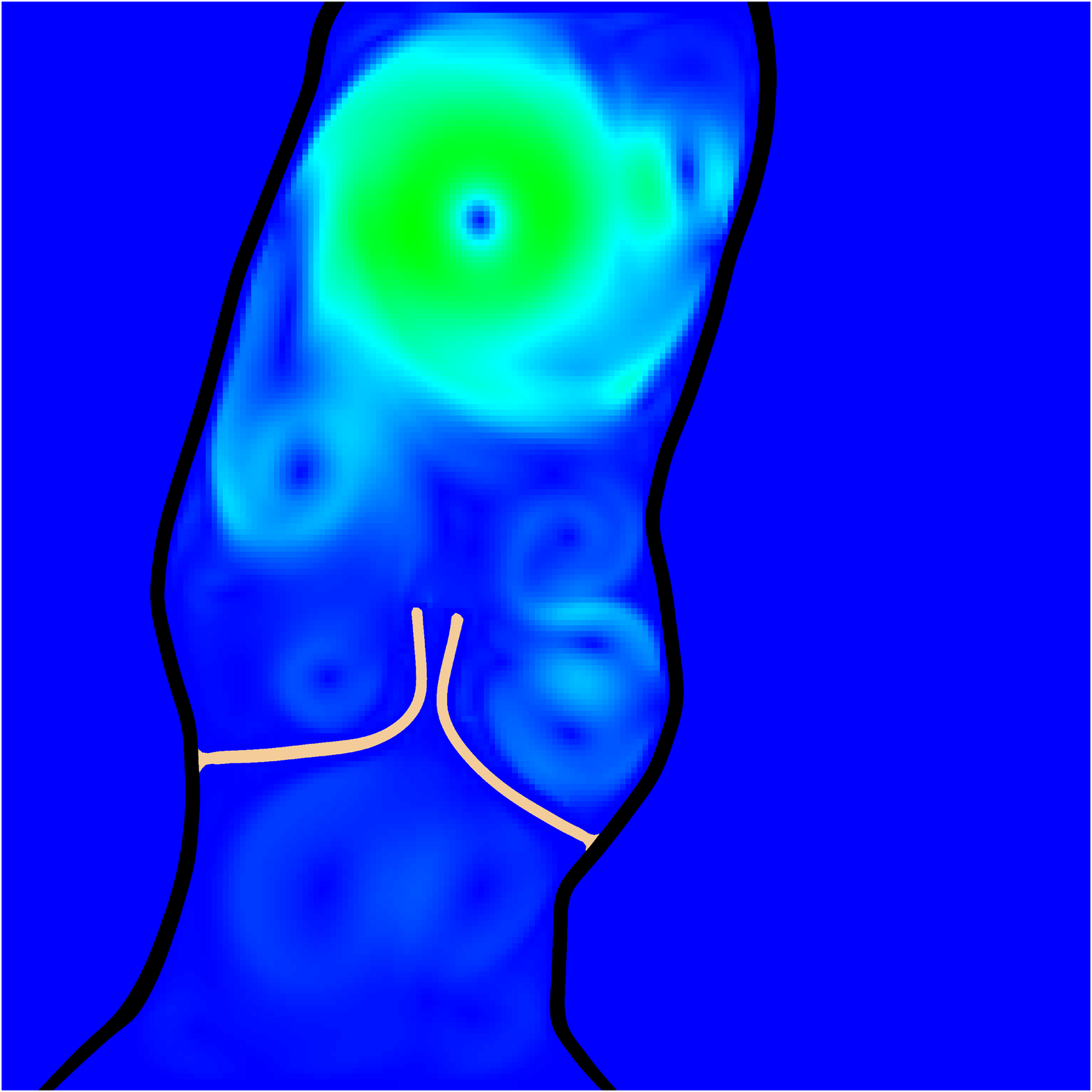} \\
\multicolumn{4}{c}{\includegraphics[width=0.5\linewidth,height=25bp]{Graphs/Results/sims/colorbarVel.eps}}
\end{tabular}
\end{center}
\caption{For a simulation in which deposition happens only on the right leaflet, snapshots during (a) middle of diastole, (b) peak systole, (c), end systole, and (d) middle of diastole the last cycle showing (top) the fluid-phase concentration $\cf$ and the surface concentration $\Cb\Js$ and (bottom) the fluid velocity magnitude. This simulation uses $\beta = 600$. Notice that the right leaflet is considerably stiffer than the left one and that the fluid concentration $\cf$ is depleted in and downstream of the aortic sinus. }\label{fig:during_sim_right}
\end{figure}

\begin{figure}
\begin{center}
\psfragscanon
\psfrag{Us}{\SI{}{\platelet\per\cm}}
\psfrag{1s}{\SI{0.220}{}}
\psfrag{2s}{\SI{0.211}{}}
\psfrag{3s}{\SI{0.201}{}}
\psfrag{4s}{\SI{0.189}{}}
\psfrag{5s}{\SI{0.173}{}}
\psfrag{6s}{\SI{0.150}{}}
\psfrag{Uf}{\SI{}{\platelet\per\cm\squared}}
\psfrag{1f}{\SI{1.000}{}}
\psfrag{2f}{\SI{0.969}{}}
\psfrag{3f}{\SI{0.933}{}}
\psfrag{4f}{\SI{0.888}{}}
\psfrag{5f}{\SI{0.831}{}}
\psfrag{6f}{\SI{0.750}{}}
\psfrag{Uv}{\SI{}{\cm\per\second}}
\psfrag{1v}{\SI{275}{}}
\psfrag{2v}{\SI{202}{}}
\psfrag{3v}{\SI{141}{}}
\psfrag{4v}{\SI{88}{}}
\psfrag{5v}{\SI{41}{}}
\psfrag{6v}{\SI{0.0}{}}
\phantomsubcaption\label{fig:during_sim_both:diastole}
\phantomsubcaption\label{fig:during_sim_both:peak_sys}
\phantomsubcaption\label{fig:during_sim_both:end_sys}
\phantomsubcaption\label{fig:during_sim_both:final}
\begin{tabular}{cccc}
\multicolumn{4}{c}{\includegraphics[width=0.5\linewidth,height=25bp,trim=0 0 0 0]{Graphs/Results/sims/colorbarRed.eps}
\includegraphics[width=0.5\linewidth,height=25bp,trim=0 0 0 0]{Graphs/Results/sims/colorbar.eps}} \\
\subfigimg[width=0.25\linewidth,hsep=0.5em,vsep=1em,pos=ul]{(a)}{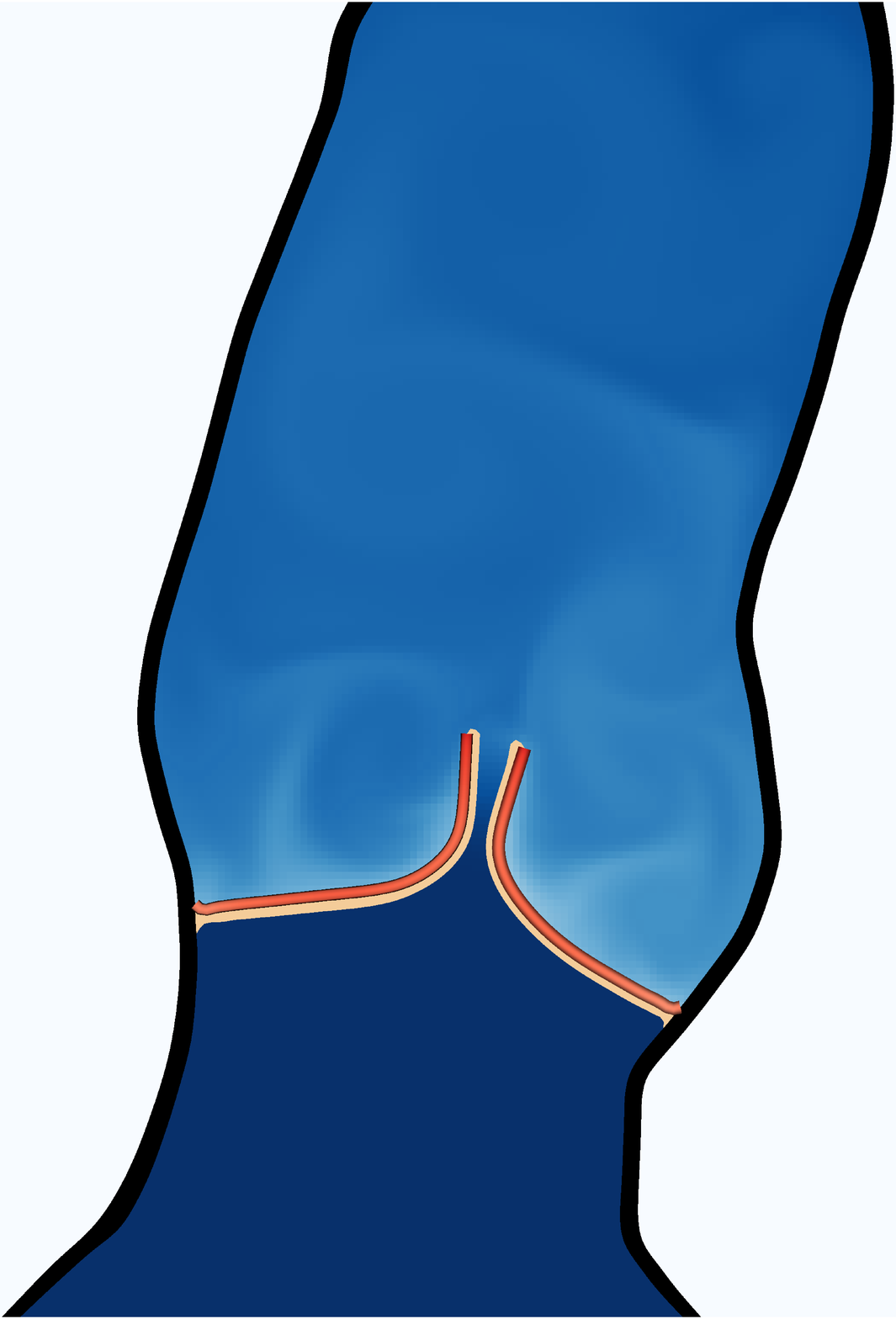}
&
\subfigimg[width=0.25\linewidth,hsep=0.5em,vsep=1em,pos=ul]{(b)}{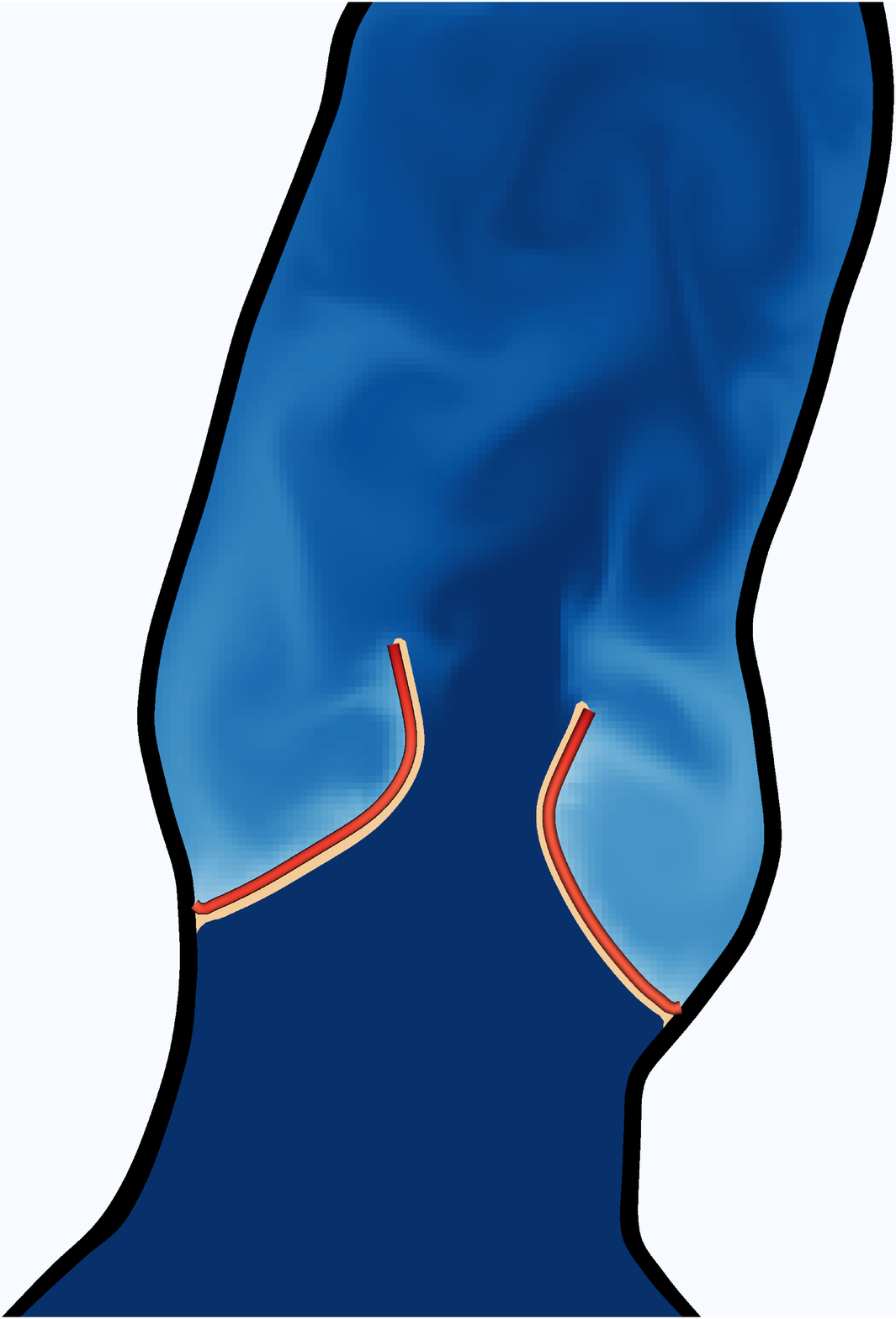}
&
\subfigimg[width=0.25\linewidth,hsep=0.5em,vsep=1em,pos=ul]{(c)}{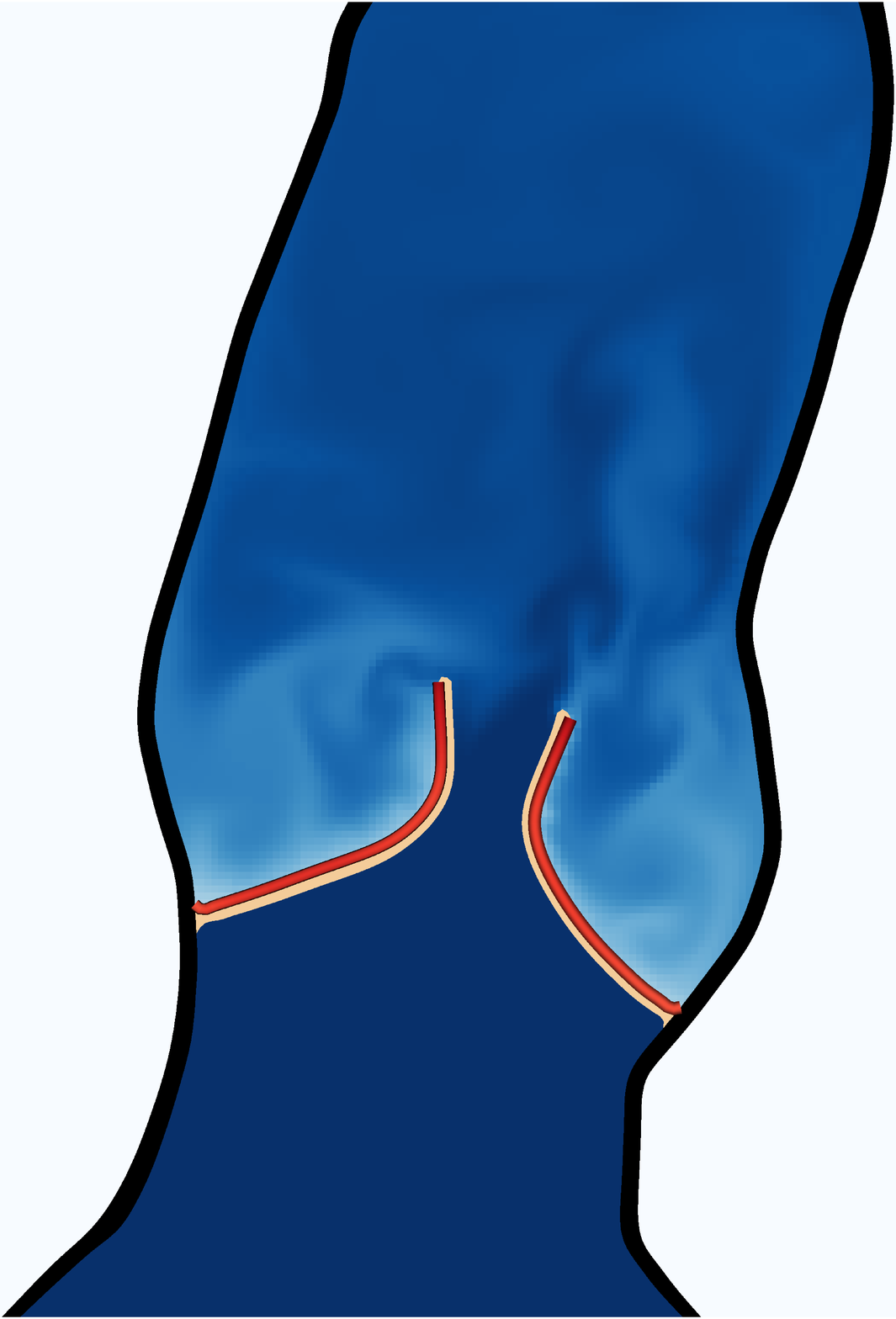}
&
\subfigimg[width=0.25\linewidth,hsep=0.5em,vsep=1em,pos=ul]{(d)}{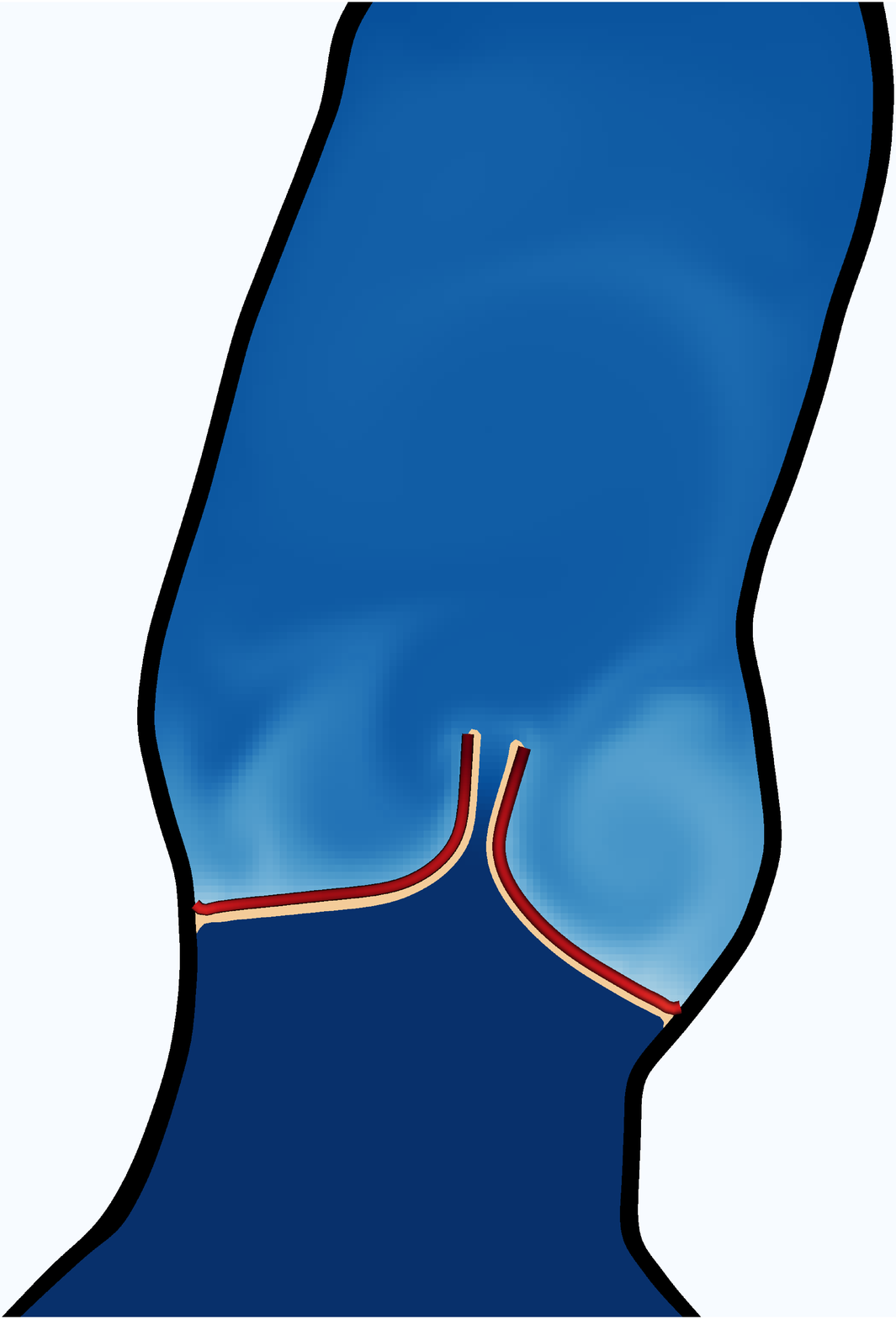} \\
\includegraphics[width=0.25\linewidth]{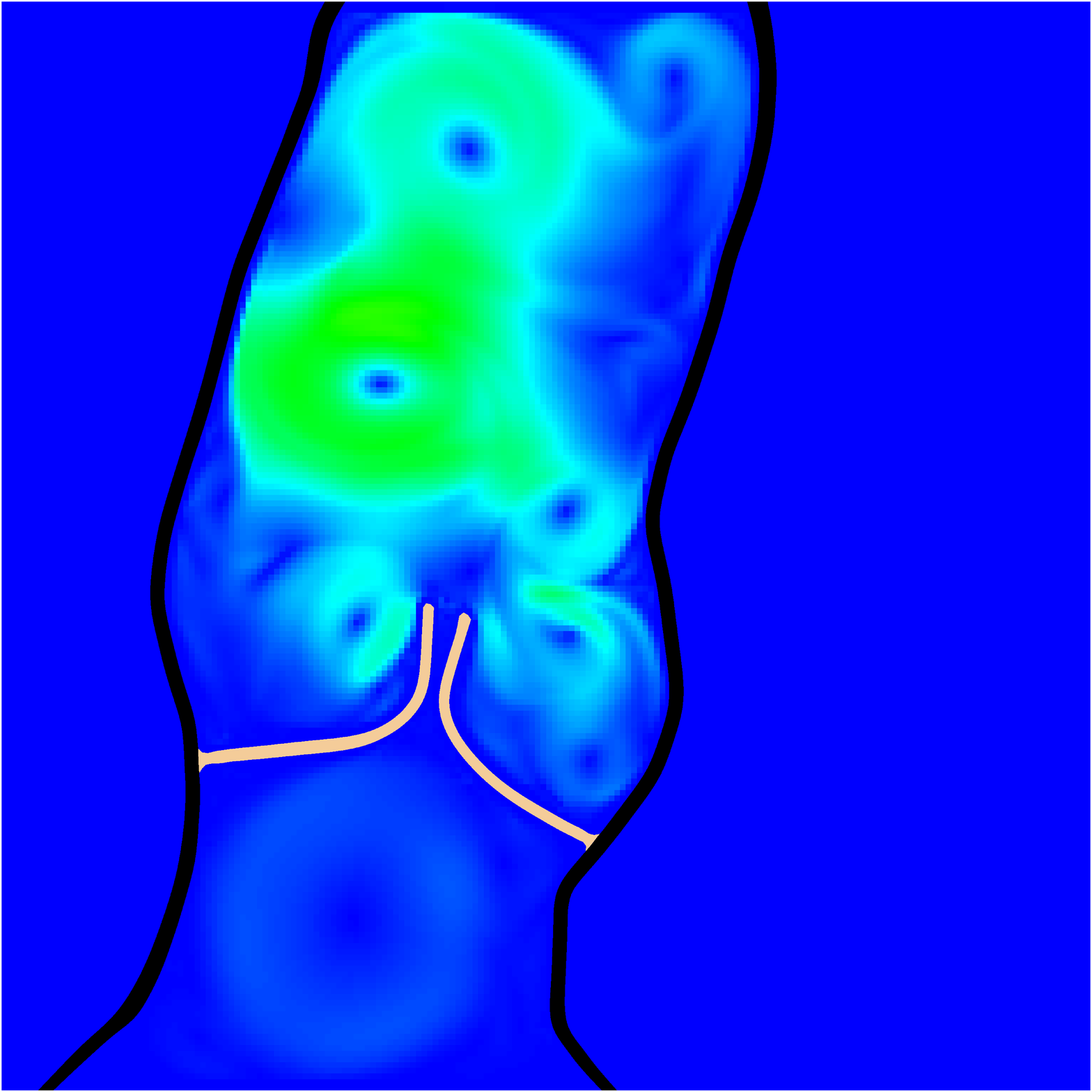}
&
\includegraphics[width=0.25\linewidth]{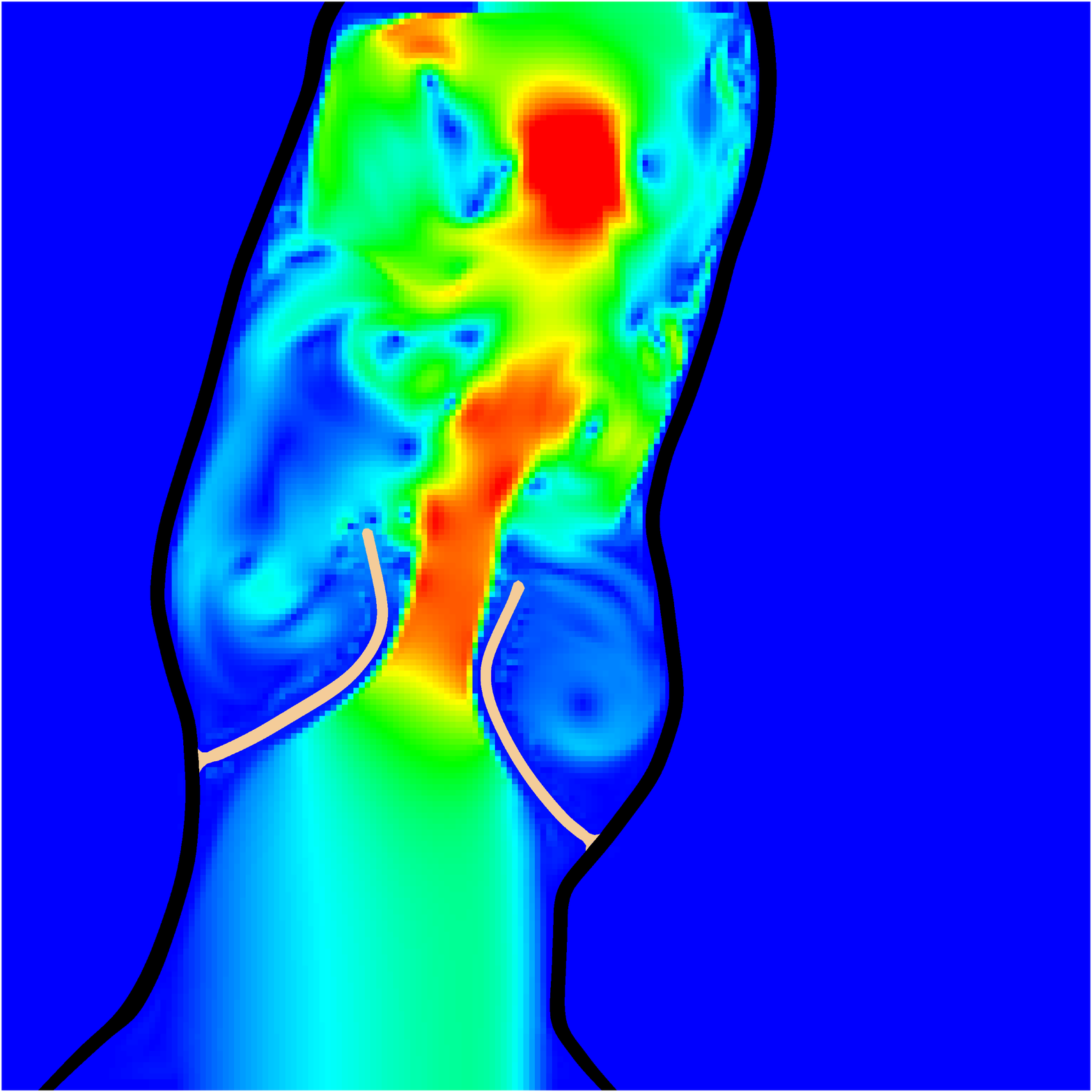}
&
\includegraphics[width=0.25\linewidth]{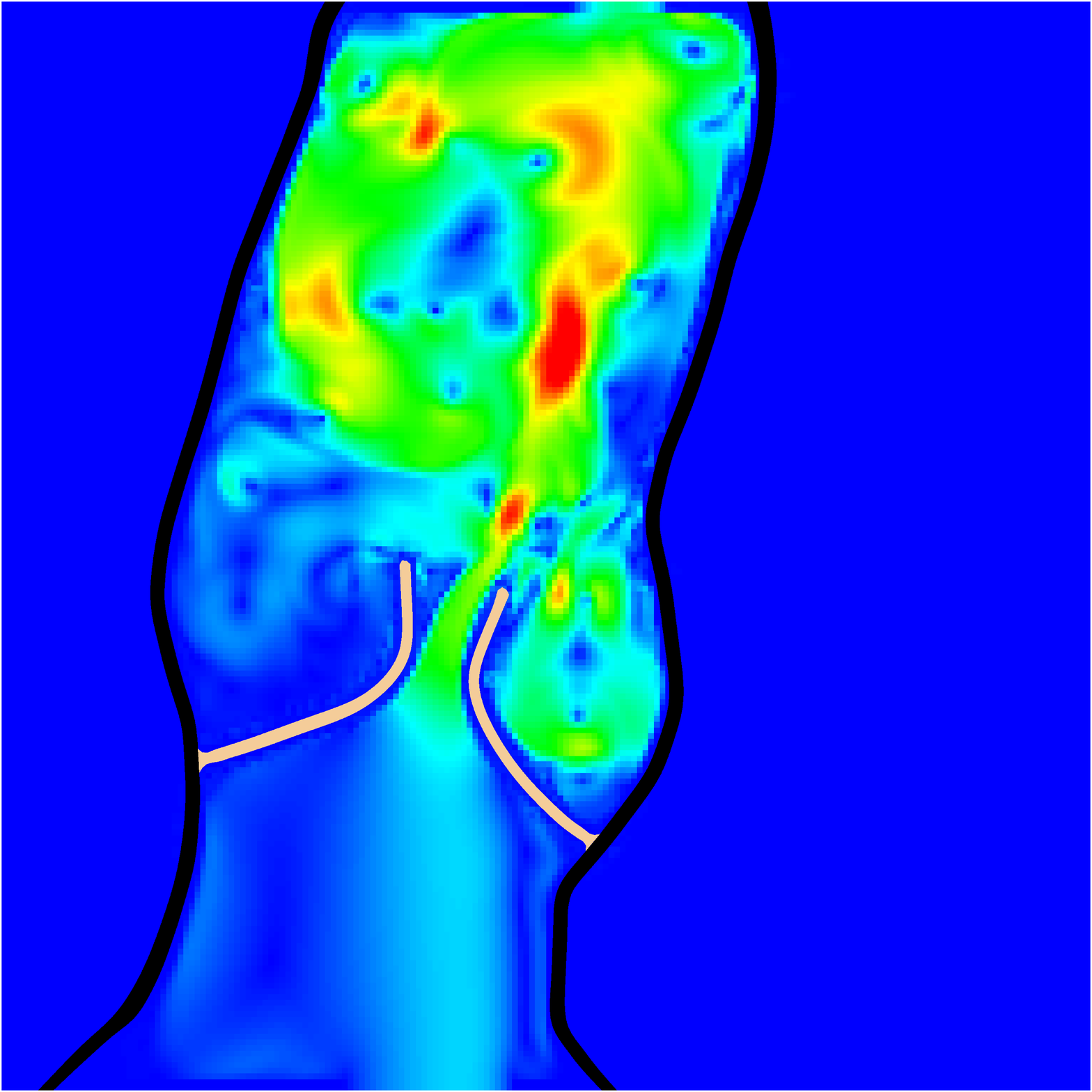}
&
\includegraphics[width=0.25\linewidth]{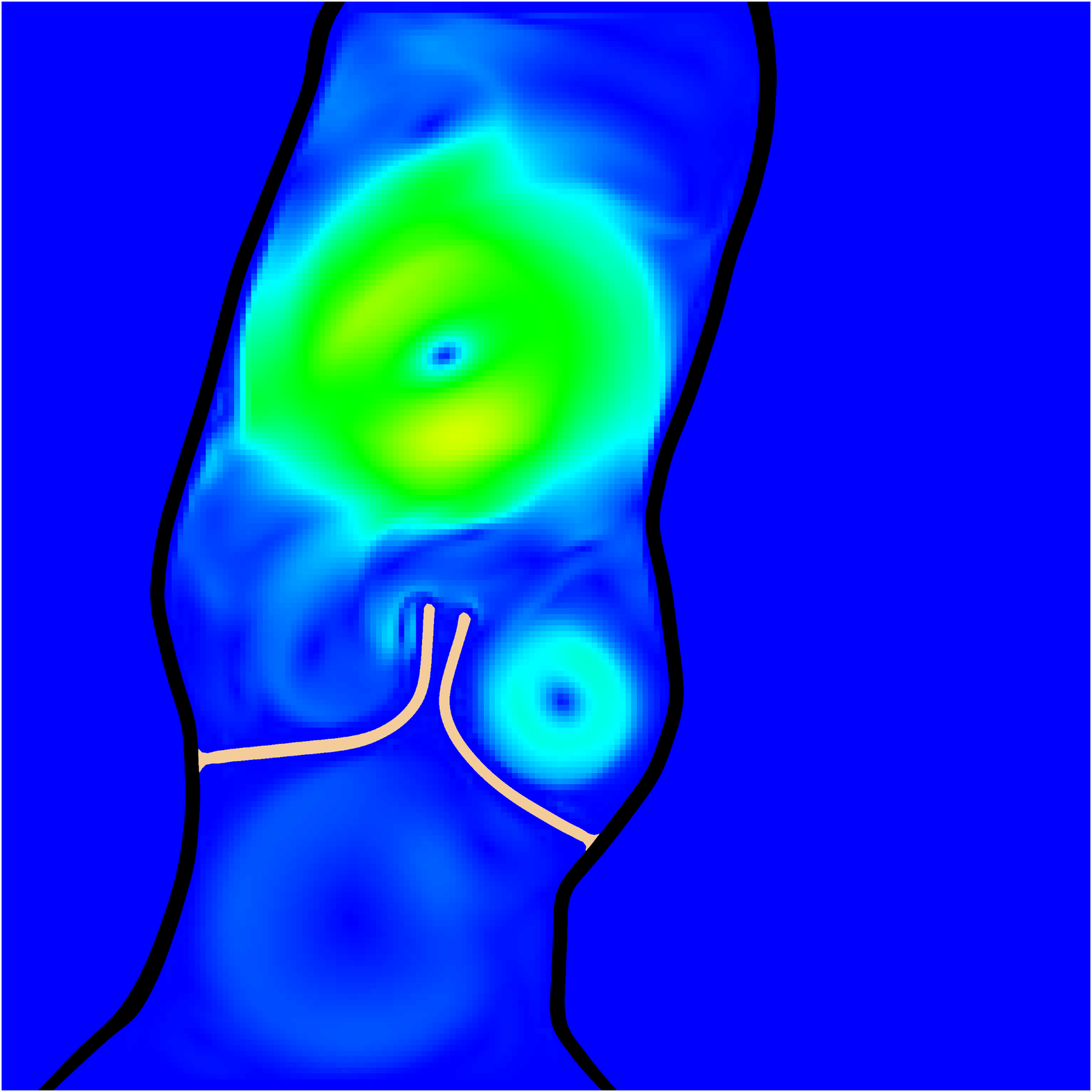} \\
\multicolumn{4}{c}{\includegraphics[width=0.5\linewidth,height=25bp]{Graphs/Results/sims/colorbarVel.eps}}
\end{tabular}
\end{center}
\caption{For a simulation in which deposition happens on both leaflets, snapshots during (a) middle of diastole, (b) peak systole, (c), end systole, and (d) middle of diastole in the last cycle showing (top) the fluid-phase concentration $\cf$ and the surface concentration $\Cb\Js$ and (bottom) the fluid velocity magnitude. This simulation uses $\beta = 600$. Notice that both leaflets stiffen in this simulation and that the peak velocity magnitudes are larger than those in \cref{fig:during_sim_right}.}\label{fig:during_sim_both}
\end{figure}

\Cref{fig:during_sim_right} shows fluid-phase concentrations $\cf$ and velocity magnitude snapshots from the last cycle of a simulation with deposition only on the right leaflet. At the times of these plots, the right leaflet is substantially stiffer than the left, and the predominant flow through the valve is shifted towards the left leaflet. \Cref{fig:during_sim_both} shows fluid-phase concentrations $\cf$ and velocity magnitude snapshots with deposition on both leaflets. Here, both leaflets become stiff and open less over time. We observe higher velocity magnitudes when deposition occurs on both leaflets as opposed to a single leaflet. 

\begin{figure}
\noindent\makebox[\textwidth][c]{
\phantomsubcaption\label{fig:area_fraction:picture}
\phantomsubcaption\label{fig:area_fraction:both}
\phantomsubcaption\label{fig:area_fraction:right_both}
\begin{minipage}{0.3\textwidth}
\centering
\subfigimg[width = 0.7\linewidth,hsep=-1em,vsep=1.5em,pos=ul]{(a)}{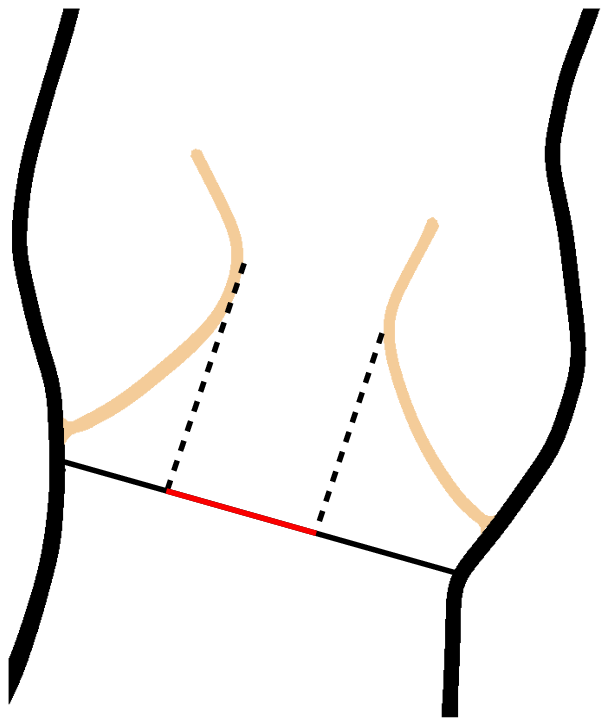}
\includegraphics[width = 0.7\linewidth]{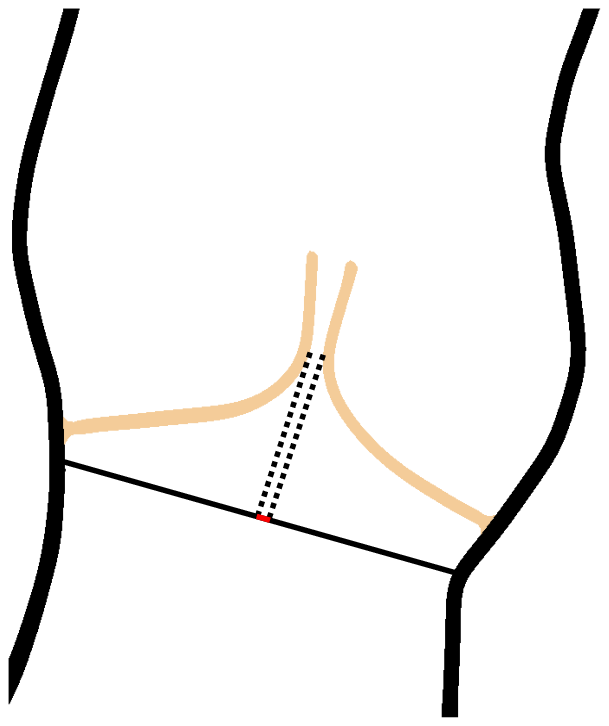}
\end{minipage}
\begin{minipage}{0.65\textwidth}
\centering
\subfigimg[width = 0.7\linewidth,hsep=-1em,vsep=1.5em,pos=ul]{(b)}{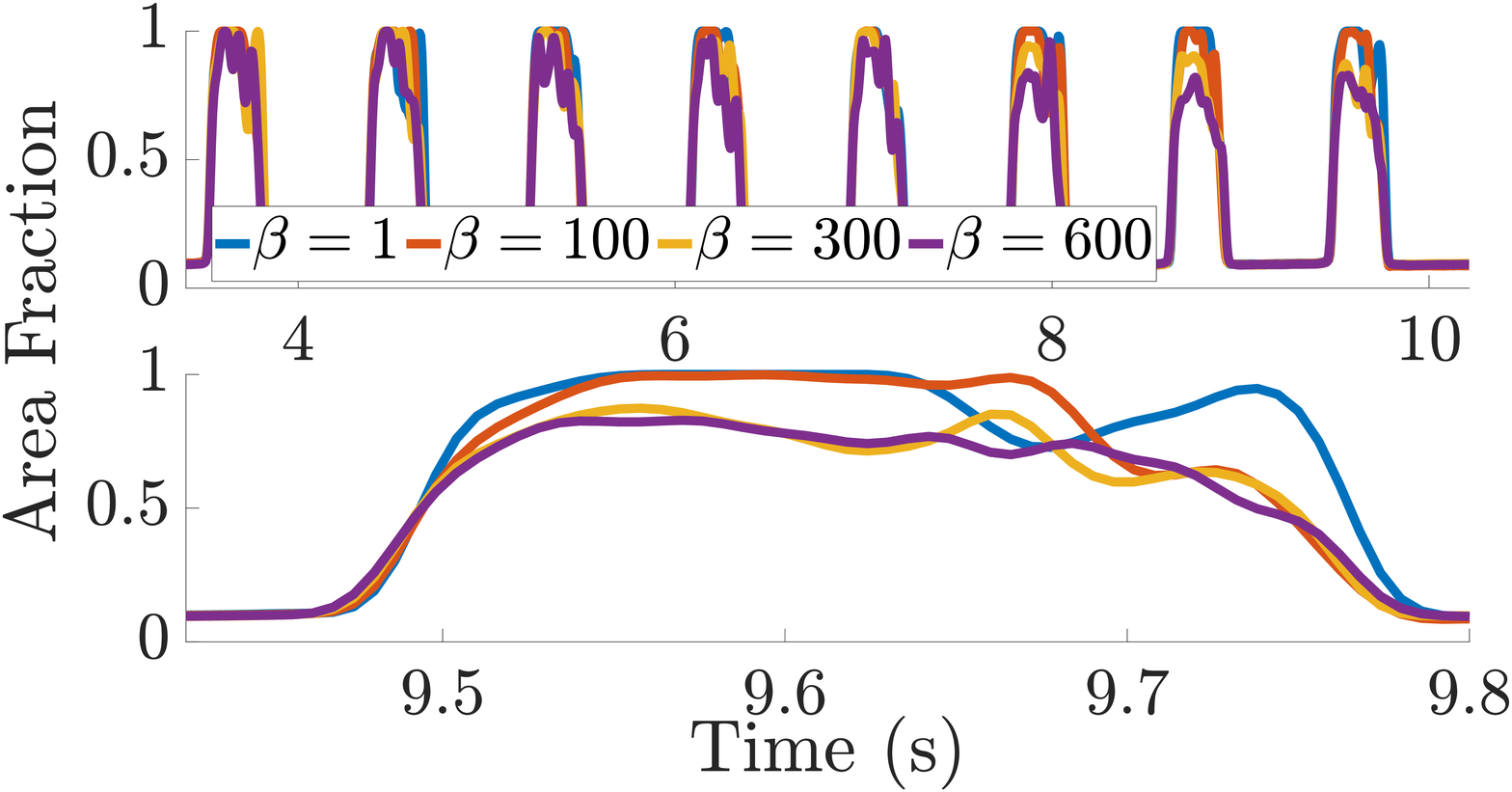}
\subfigimg[width = 0.7\linewidth,hsep=-1em,vsep=1.5em,pos=ul]{(c)}{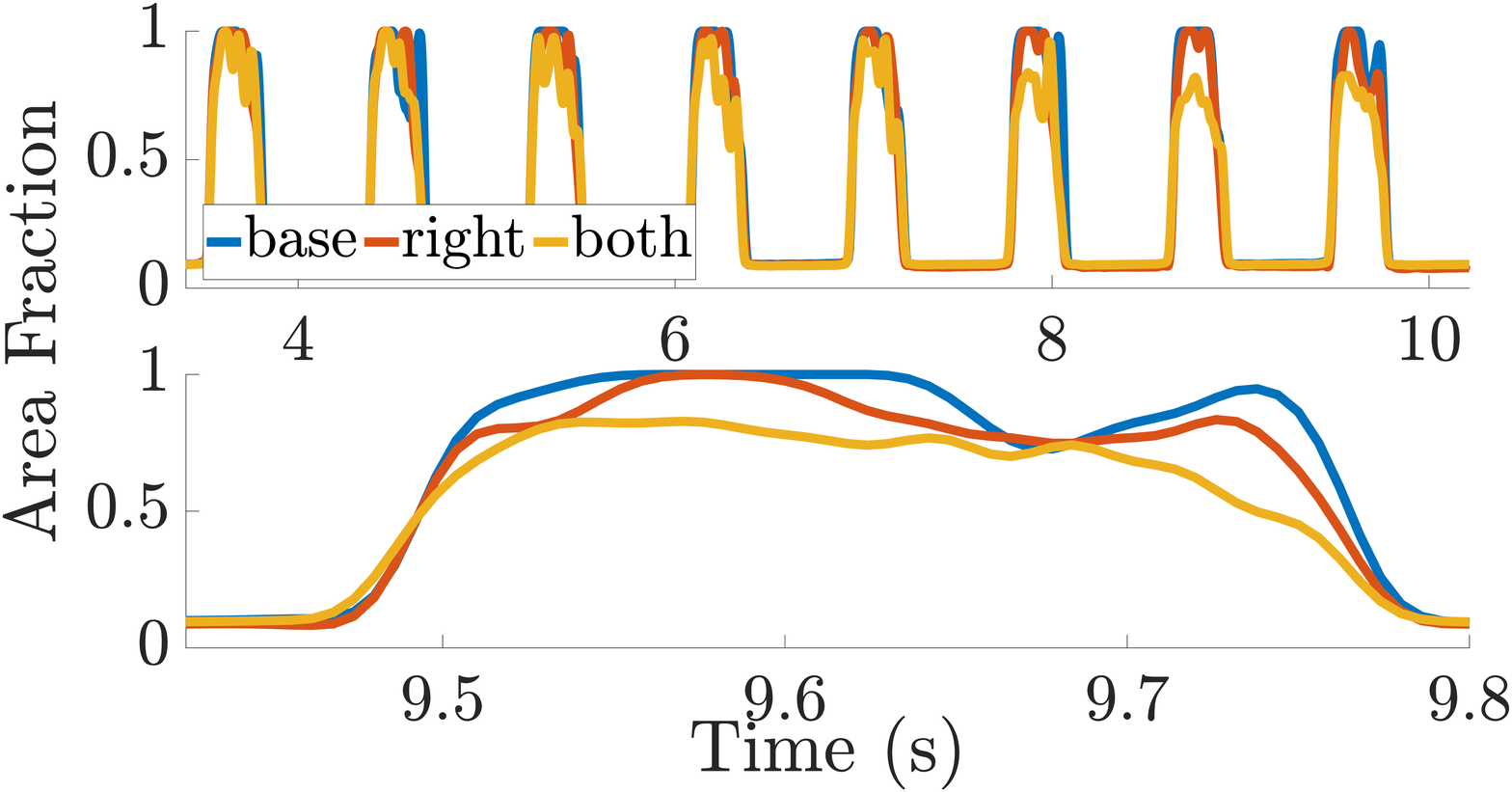}
\end{minipage}
}
\caption{The computed normalized open valve area over time. Panel (a) depicts the computation of the normalized open valve area. The open area of the valve is projected onto the valve ring, and then normalized by the area of the fully opened valve. Panel (b) depicts the normalized open valve area over time for accumulation on both leaflets as $\beta$ increases during the entire simulation (top) and during only the last cycle (bottom). Notice that as the leaflets get stiffer, the normalized open valve area decreases. Panel (c) depicts the normalized open valve area over time for accumulation on the right leaflet or both leaflets during all the cycles (top) and during only the last cycle (bottom). Notice that if accumulation occurs on both leaflets, the normalized open valve area decreases more than if accumulation occurs on a single leaflet.} \label{fig:area_fraction}
\end{figure}

To assess the opening of the valve, we project the leaflets onto the valve ring, as shown in \cref{fig:area_fraction:picture}. The opening area is then normalized by the area of the fully open valve. \Cref{fig:area_fraction:both} shows the normalized open valve area over each cycle for deposition on both leaflets, as we increase the maximum stiffness factor $\beta$. For lower maximal stiffnesses, we observe similar normalized open valve areas compared to a simulation with no deposition. For larger maximum stiffness, we observe a smaller normalized open valve area as more deposition occurs. \Cref{fig:area_fraction:right_both} compares the normalized open valve area for deposition on both leaflets versus on only the right leaflet for the same maximum stiffness. When deposition occurs on only the right leaflet, the normalized open valve area still decreases compared to no accumulation, but does not realize as dramatic reductions as when deposition on both leaflets is allowed. The left leaflet, which has a constant stiffness over time, compensates and opens more as the right leaflet stiffens.

\begin{figure}
\begin{center}
\psfragscanon
\phantomsubcaption{\label{fig:concentration:accumulation}}
\phantomsubcaption\label{fig:concentration:both}
\phantomsubcaption\label{fig:concentration:both_close}
\phantomsubcaption\label{fig:concentration:right_both_close}
\psfrag{n1}{0.200}
\psfrag{n2}{0.205}
\psfrag{n3}{0.210}
\psfrag{n4}{0.215}
\psfrag{n5}{0.220}
\begin{minipage}{0.44\textwidth}
\centering
\hspace{-2em}\subfigimg[width = 0.65\linewidth,hsep=-1.25em,vsep=1em,pos=ul]{(a)}{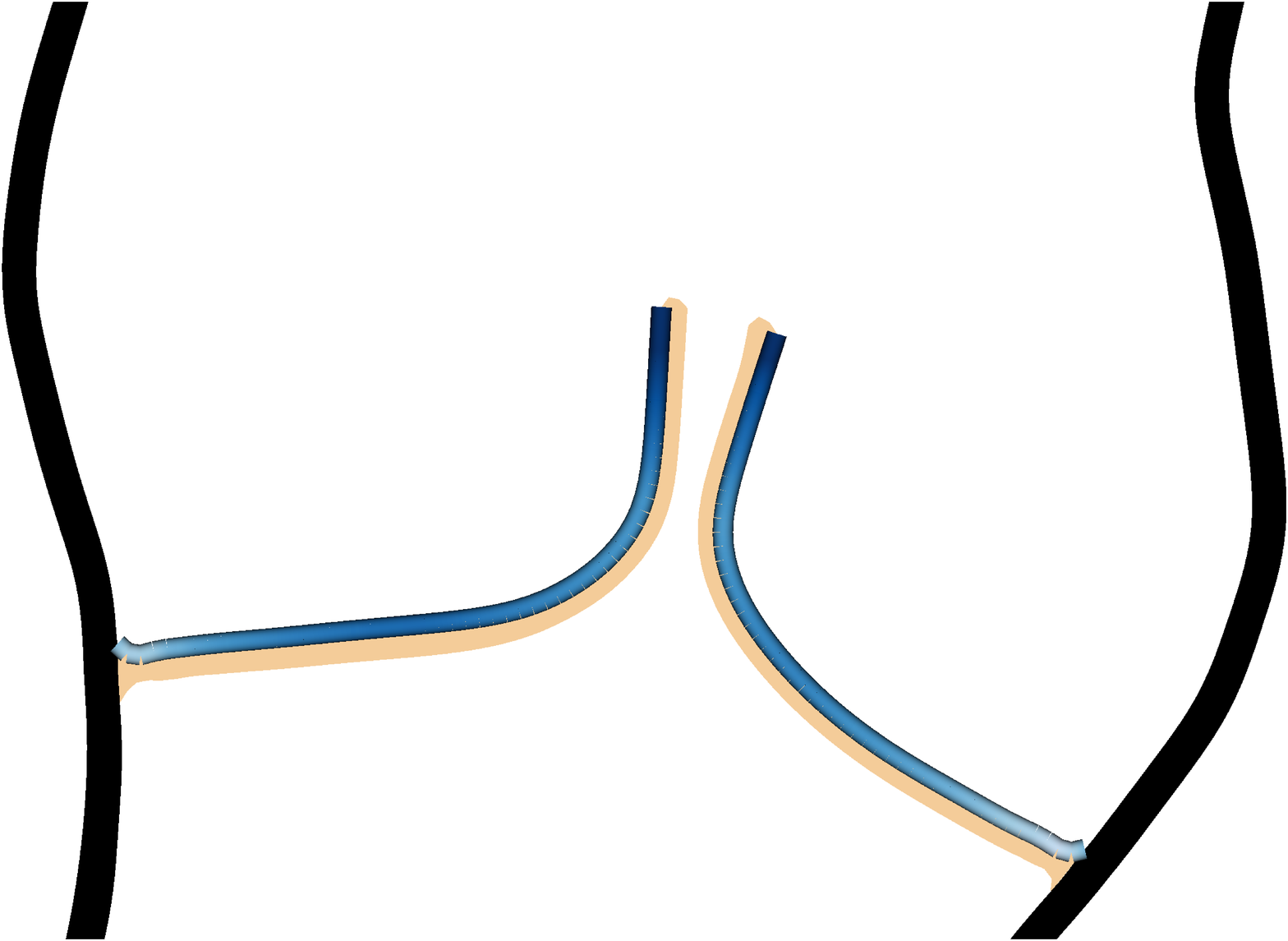}
\includegraphics[width=0.1\linewidth,height=93bp]{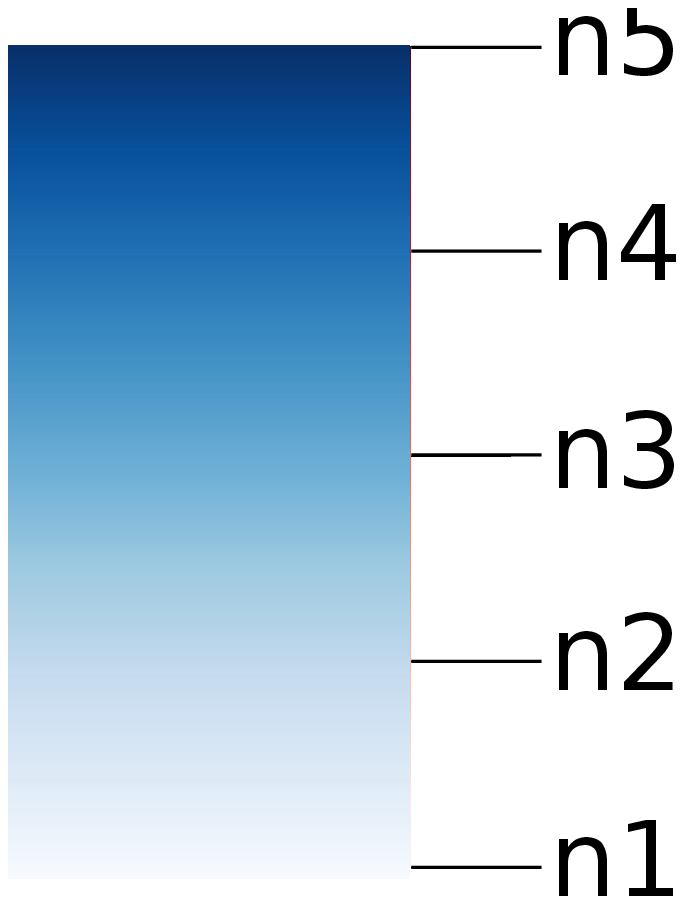}
\subfigimg[width = 0.9\linewidth,hsep=-1em,vsep=1.5em,pos=ul]{(c)}{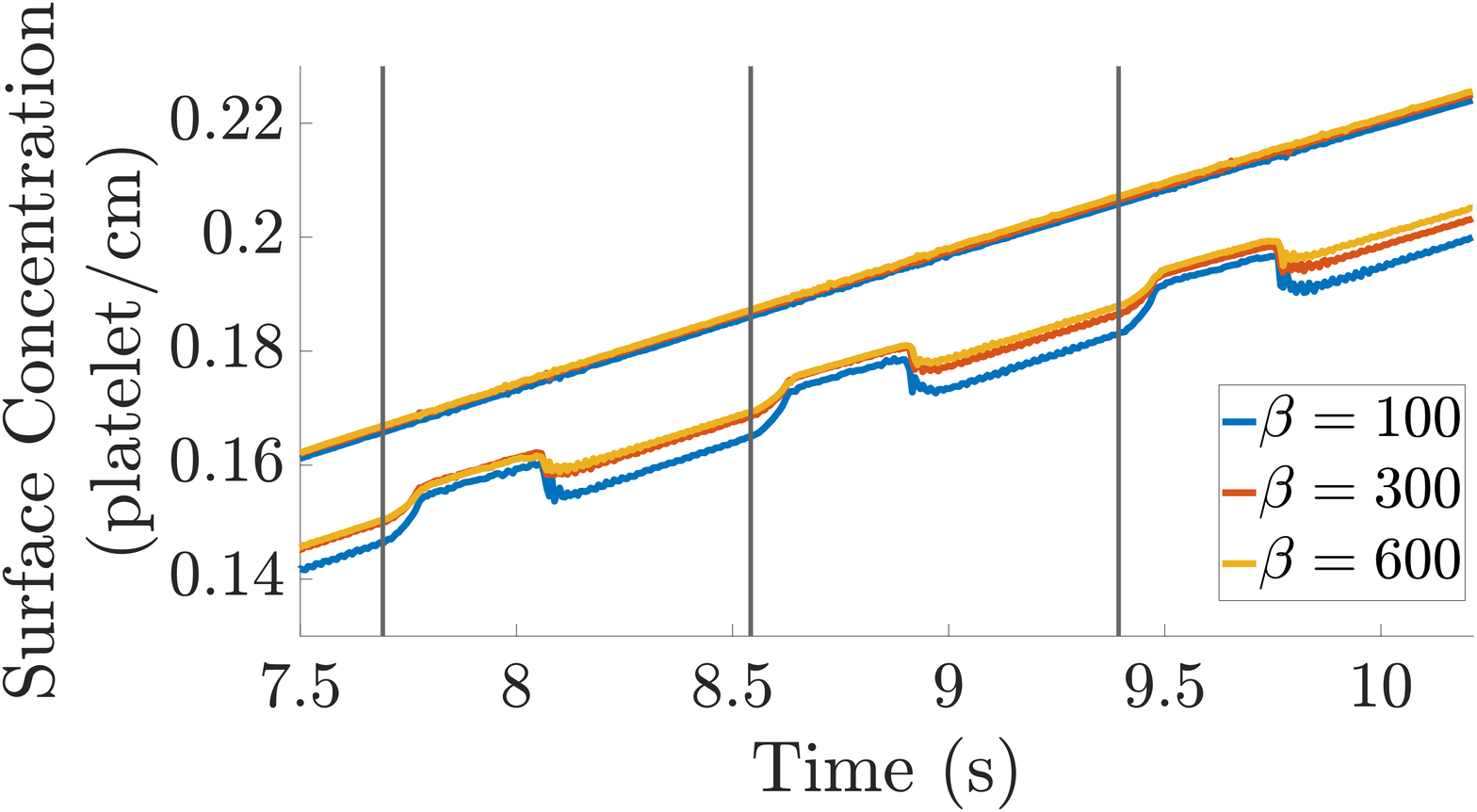}

\end{minipage}
\begin{minipage}{0.44\textwidth}
\subfigimg[width = 0.9\linewidth,hsep=-1em,vsep=1.5em,pos=ul]{(b)}{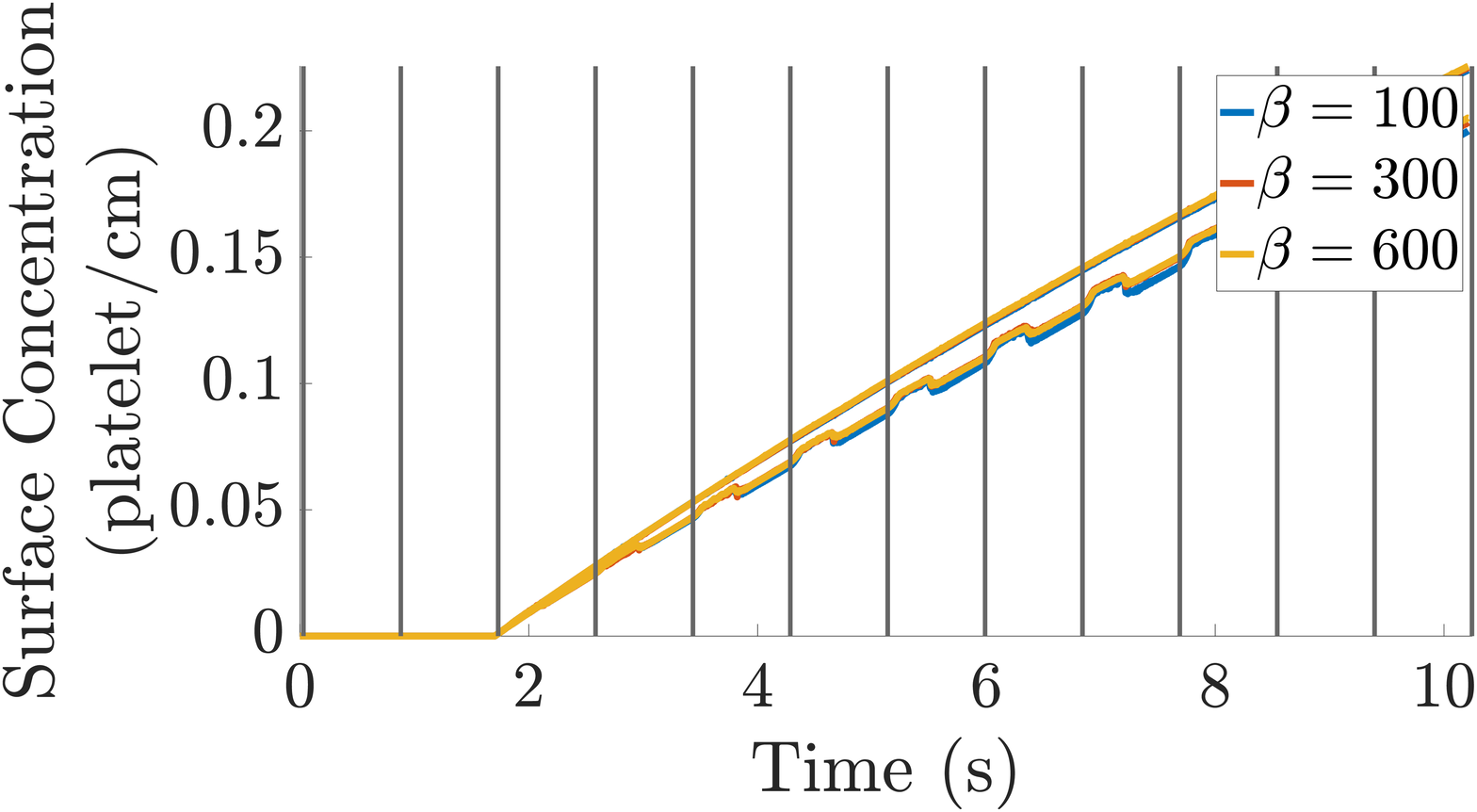}
\subfigimg[width = 0.9\linewidth,hsep=-1em,vsep=1.5em,pos=ul]{(d)}{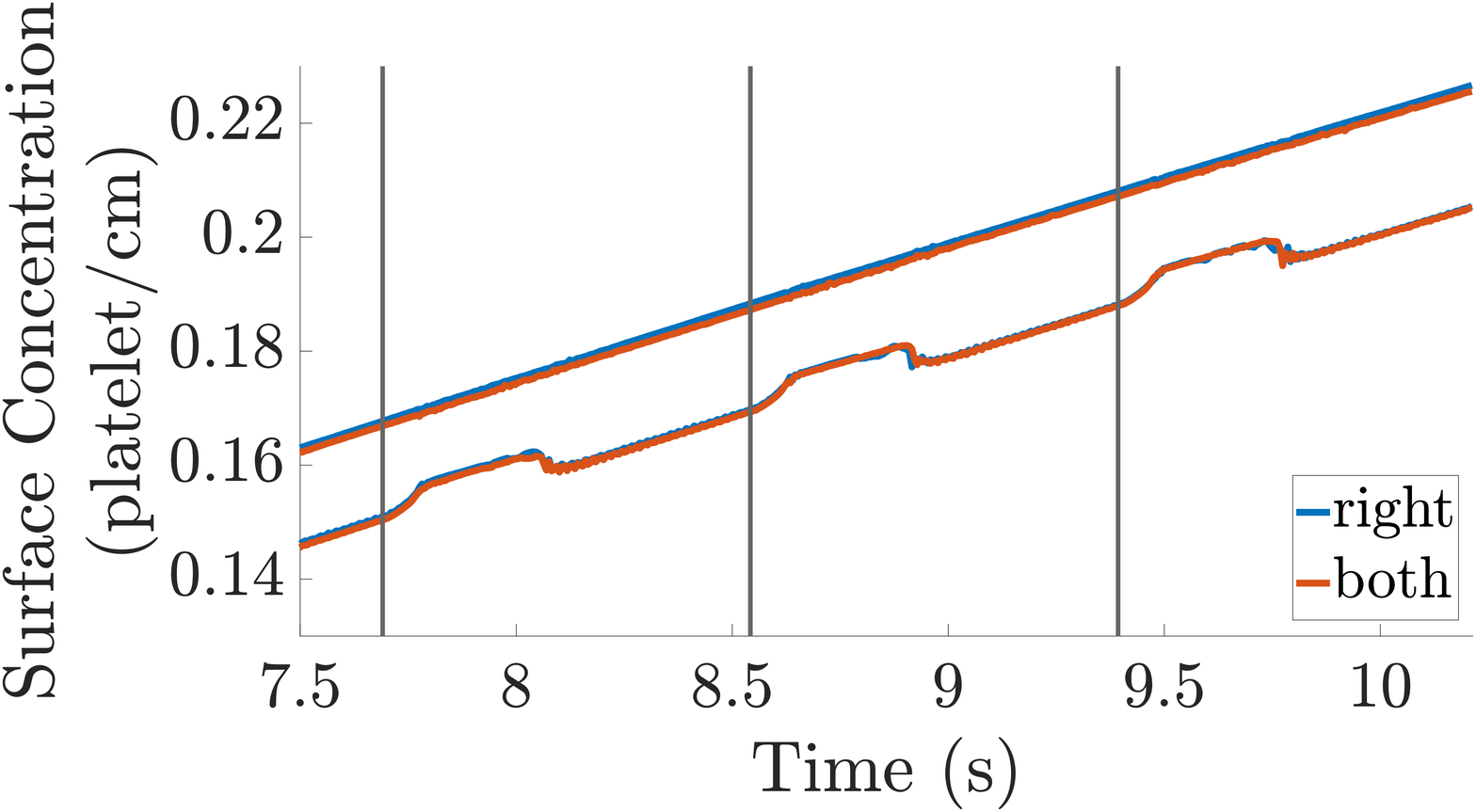}
\end{minipage}
\end{center}
\caption{(a) The surface concentration $\Cb\Js$ along the leafets at the end of the simulation with $\beta = 600$. (b) The minimum and maximum surface concentration $\Cb\Js$ on the leaflets. Panel (c) highlights the accumulation over the last three cycles for accumulation on both leaflets for three values of $\beta$. In panel (d), the accumulation is shown for accumulation only on the right leaflet versus both leaflets with $\beta = 600$. The vertical lines denote the beginning of systole. There is a consistent accumulation of material on the leaflets. The jumps in the minimum concentration are due to the the change in $\Js$, which changes most where the leaflets attach to the aortic wall.}
\label{fig:concentration}
\end{figure}

\Cref{fig:concentration} shows the maximum and minimum accumulations $\Cb \Js$ across the leaflet. Because the diffusion coefficient is large compared to the reaction rates, there is always sufficient fluid-phase platelets to bind to the leaflets, and accordingly, a consistent rate of binding to the surface and an increasing surface concentration as the simulation progresses. By the end of the simulation, bound platelets occupy approximately 23\% of the carrying capacity. The minimum concentration periodically jumps while the maximum concentration is monotonically increasing. The periodic jumps are due to the physical location of the minimum and maximum which affects the amount of platelets per current area. The minimum occurs near the position where the leaflet attaches to the aortic wall. This location sees the largest changes in area as the valve opens. The maximum concentration is found on the tips of the leaflet, which move through regions of high fluid-phase concentration. The tip of the leaflets deform less than the rest of the leaflet upon opening and closing of the valve, leading to a steadily increasing surface bound concentration field. While the fluid-phase concentration $\cf$ is not completely depleted, we do observe reductions in the concentration than initially.

\begin{figure}
\begin{center}
\phantomsubcaption\label{fig:fluid_flow:base}
\phantomsubcaption\label{fig:fluid_flow:100x}
\phantomsubcaption\label{fig:fluid_flow:300x}
\phantomsubcaption\label{fig:fluid_flow:600x}
\phantomsubcaption\label{fig:fluid_flow:right}
\psfragscanon
\psfrag{1}{\SI{50}{}}
\psfrag{2}{\SI{40}{}}
\psfrag{3}{\SI{30}{}}
\psfrag{4}{\SI{20}{}}
\psfrag{5}{\SI{10}{}}
\psfrag{6}{\SI{0}{}}
\psfrag{sc}{\SI{}{\cm\per\second}}
\hspace{-3em}
\subfigimg[width = 0.38\linewidth,hsep=1.75em,vsep=1.5em,pos=ul]{(a)}{}
\hspace{-2.5em}
\subfigimg[width = 0.38\linewidth,hsep=1.75em,vsep=1.5em,pos=ul]{(b)}{}
\hspace{-2.5em}
\subfigimg[width = 0.38\linewidth,hsep=1.75em,vsep=1.5em,pos=ul]{(c)}{} \\
\hspace{-3em}
\subfigimg[width = 0.38\linewidth,hsep=1.75em,vsep=1.5em,pos=ul]{(d)}{}
\hspace{-2.5em}
\subfigimg[width = 0.38\linewidth,hsep=1.75em,vsep=1.5em,pos=ul]{(e)}{}
\includegraphics[width=0.08\linewidth,height=109bp,trim=0 -45 0 90]{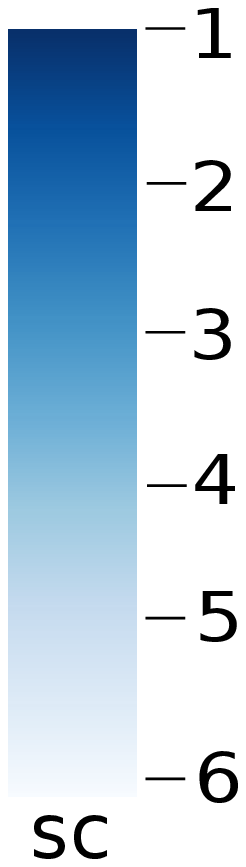}
\end{center}
\caption{The magnitude of the velocity field at peak systole during the last cycle for (a) no accumulation, (b) $\beta = 100$ with accumulation on both leaflets, (c) $\beta = 300$ with accumulation on both leaflets, (d) $\beta = 600$ with accumulation on both leaflets, and (e) $\beta = 600$ with accumulation on only the right leaflet. Notice the presence of a vortex in the right sinus that is absent if accumulation occurs only on the right leaflet. The choice of colorbar is intended to highlight vortex formation in the sinus region. We observe peak flow rates through the valve of \SI{170}{\cm\per\second} for $\beta = 1$ and \SI{275}{\cm\per\second} for $\beta = 600$.}\label{fig:fluid_flow}
\end{figure}

\Cref{fig:fluid_flow} shows the velocity magnitude during the final cycle near peak systole. We observe a vortex in the sinus region that grows in strength as we increase the maximum stiffness. This vortex is not present when accumulation occurs exclusively on the right leaflet, but is present when there is no accumulation.

\subsection{Pressures and Flow Rates}\label{sec:RESULTS:stroke_vol}

\begin{figure}
\begin{center}
\phantomsubcaption\label{fig:pressures:both}
\phantomsubcaption\label{fig:pressures:right_both}
\psfragscanon
\psfrag{LVOT - BE}{\scalebox{0.45}{LVOT: $\beta = 1$}}
\psfrag{aorta - BE}{\scalebox{0.45}{aorta: $\beta = 1$}}
\psfrag{LVOT - BETA}{\scalebox{0.45}{LVOT: $\beta = 600$}}
\psfrag{aorta - BETA}{\scalebox{0.45}{aorta: $\beta = 600$}}
\subfigimg[width=0.49\linewidth,hsep=-0.75em,vsep=2em,pos=ul]{(a)}{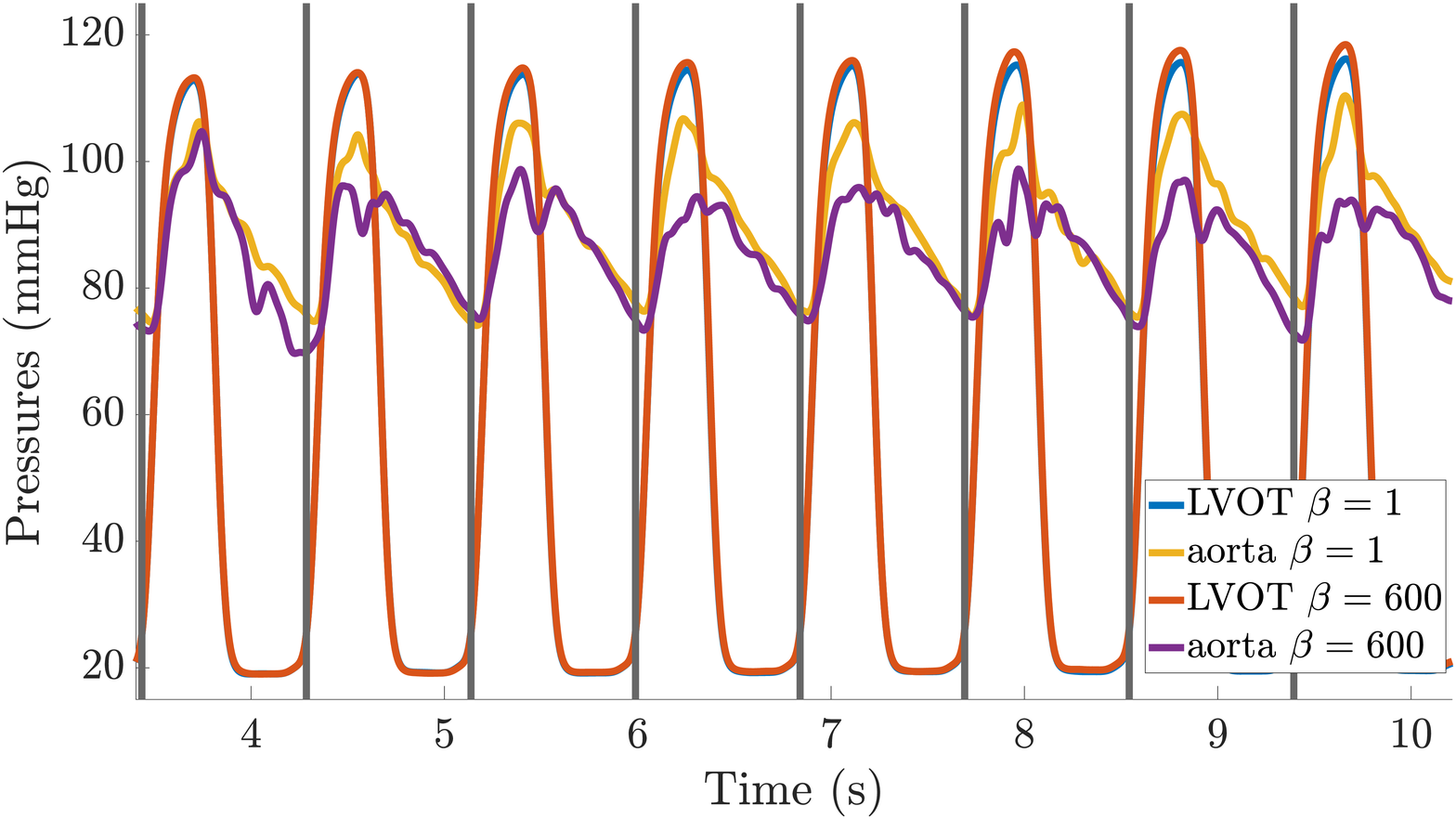}
\subfigimg[width=0.49\linewidth,hsep=-0.75em,vsep=2em,pos=ul]{(b)}{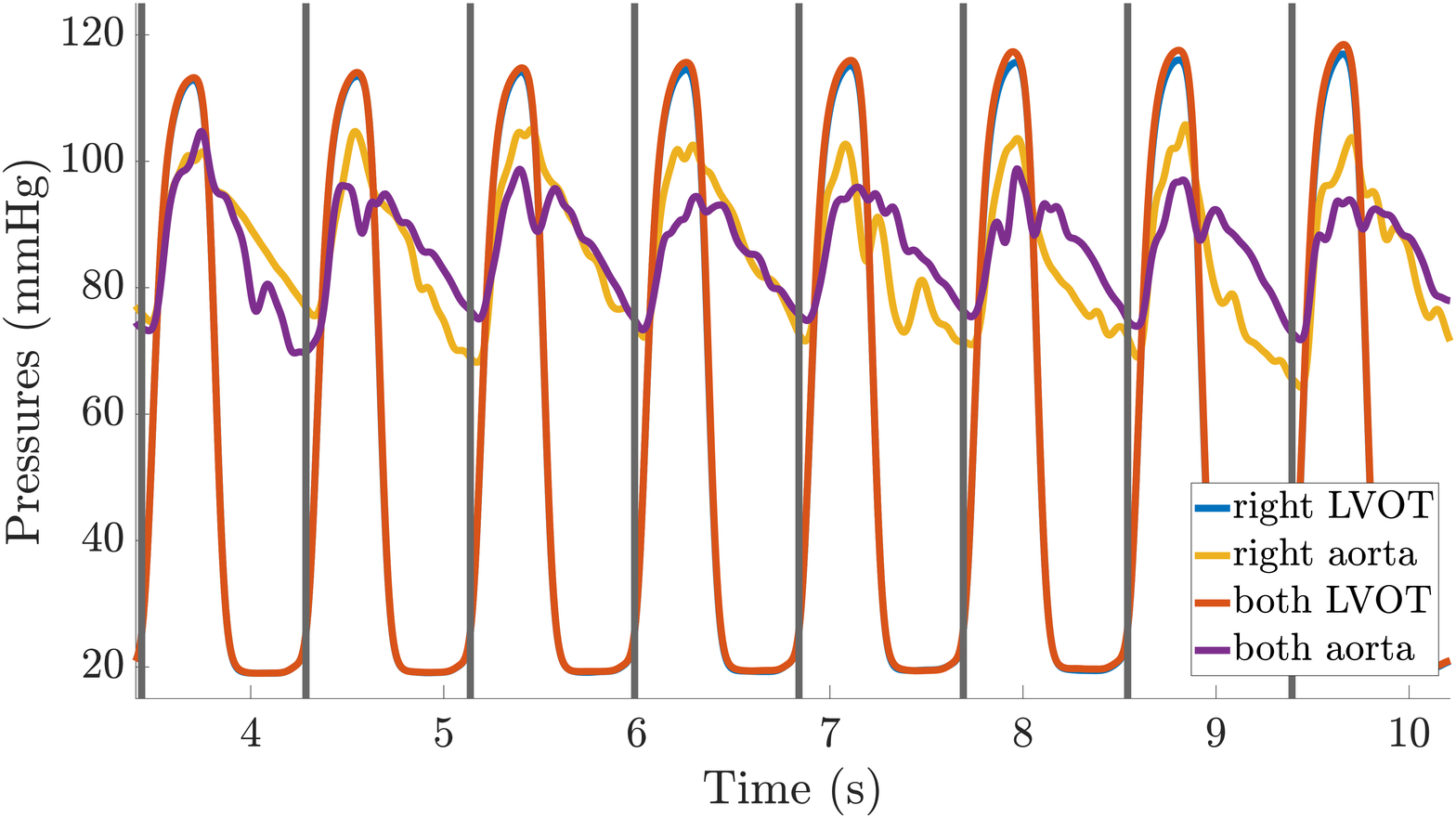}
\end{center}
\caption{The pressures just upstream and downstream of the valve. Panel (a) compares the pressures for no accumulation and accumulation on both leaflets with $\beta = 600$. Panel (b) compares pressures for accumulation on only the right leaflet and both leaflets, with $\beta = 600$. While the upstream pressure increases mildly, the downstream pressures decrease by between $5$ and $\SI{10}{\mmHg}$.}\label{fig:pressures}
\end{figure}

\begin{figure}
\begin{center}
\phantomsubcaption\label{fig:pressures_close:aorta}
\phantomsubcaption\label{fig:pressures_close:LVOT}
\subfigimg[width=0.49\linewidth,hsep=-0.75em,vsep=2em,pos=ul]{(a)}{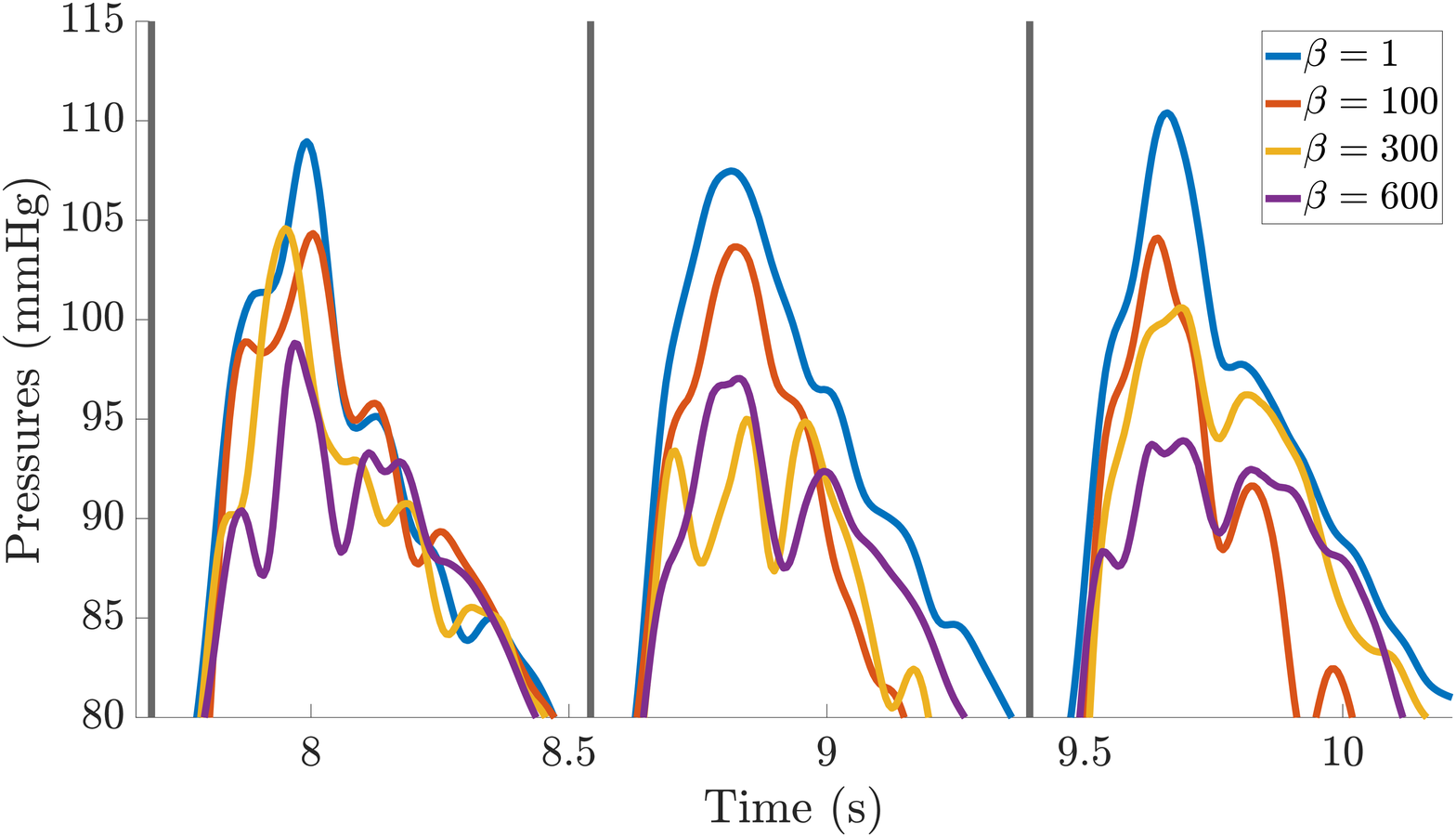}
\subfigimg[width=0.49\linewidth,hsep=-0.75em,vsep=2em,pos=ul]{(b)}{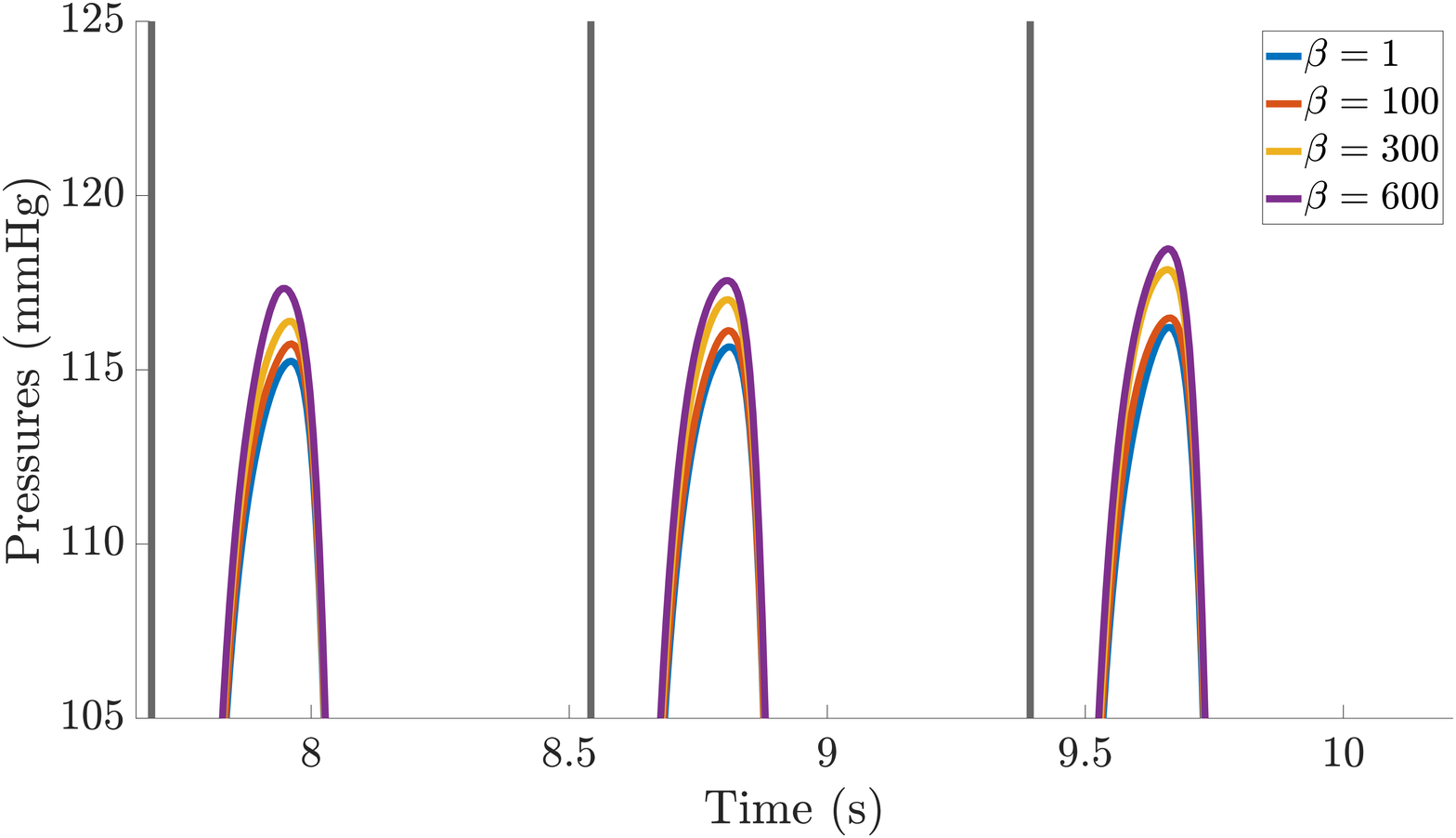}
\end{center}
\caption{The pressures during the last three cycles for accumulation on both leaflets with $\beta = 1$, 100, 300, and 600. Note that $\beta = 1$ corresponds to a baseline where the leaflet stiffness is constant. The aortic side pressures are shown in panel (a). The left ventricle side pressures are shown in panel (b). While the upstream pressure increases by less than $\SI{5}{\mmHg}$, the downstream pressures decrease by between $5$ and $\SI{15}{\mmHg}$.} \label{fig:pressures_close}
\end{figure}

Here we quantify the valve's resistance to the flow at different maximum stiffnesses. We measure the pressures at locations just upstream and downstream of the valve. We use a Gaussian filter to smooth the curves in both space and time, yielding the results shown in \cref{fig:pressures,fig:pressures_close}. We observe marginal increases in the pressures upstream of the valve of about $\SIrange{1}{4}{\mmHg}$; however, there are sharp decreases of about $\SIrange{5}{15}{\mmHg}$ in the aortic pressures downstream of the leaflets as we increase stiffness. A similar trend is observed with deposition on only the right leaflet, although the differences are not as pronounced.

\begin{table}
\centering
\caption{Pressures just upstream of the leaflets (LVOT) and downstream of the leaflets (aorta) during peak systole. While we observe increases of \SIrange{1}{2}{\mmHg} in the pressure upstream of the valve, there are greater decreases of \SIrange{10}{15}{\mmHg} in the pressure downstream of the valve.}
    \begin{tabular}{c | c | c c c | c}
      & Base & \multicolumn{3}{c|}{Both Leaflets} & Right Leaflet \\
      &      & $\beta = 100$ & $\beta = 300$ & $\beta = 600$ & $\beta = 600$ \\
      \hline
LVOT  & \SI{117}{\mmHg}  & \SI{117}{\mmHg} & \SI{118}{\mmHg} & \SI{119}{\mmHg}  & \SI{117}{\mmHg}  \\
Aorta & \SI{111}{\mmHg}  & \SI{105}{\mmHg} & \SI{101}{\mmHg} & \SI{94}{\mmHg} & \SI{104}{\mmHg}
    \end{tabular}
    \label{results:tab:pressures}
\end{table}

\begin{figure}
\subfigimg[width=0.5\columnwidth,hsep=-0.75em,vsep=2em,pos=ul]{(a)}{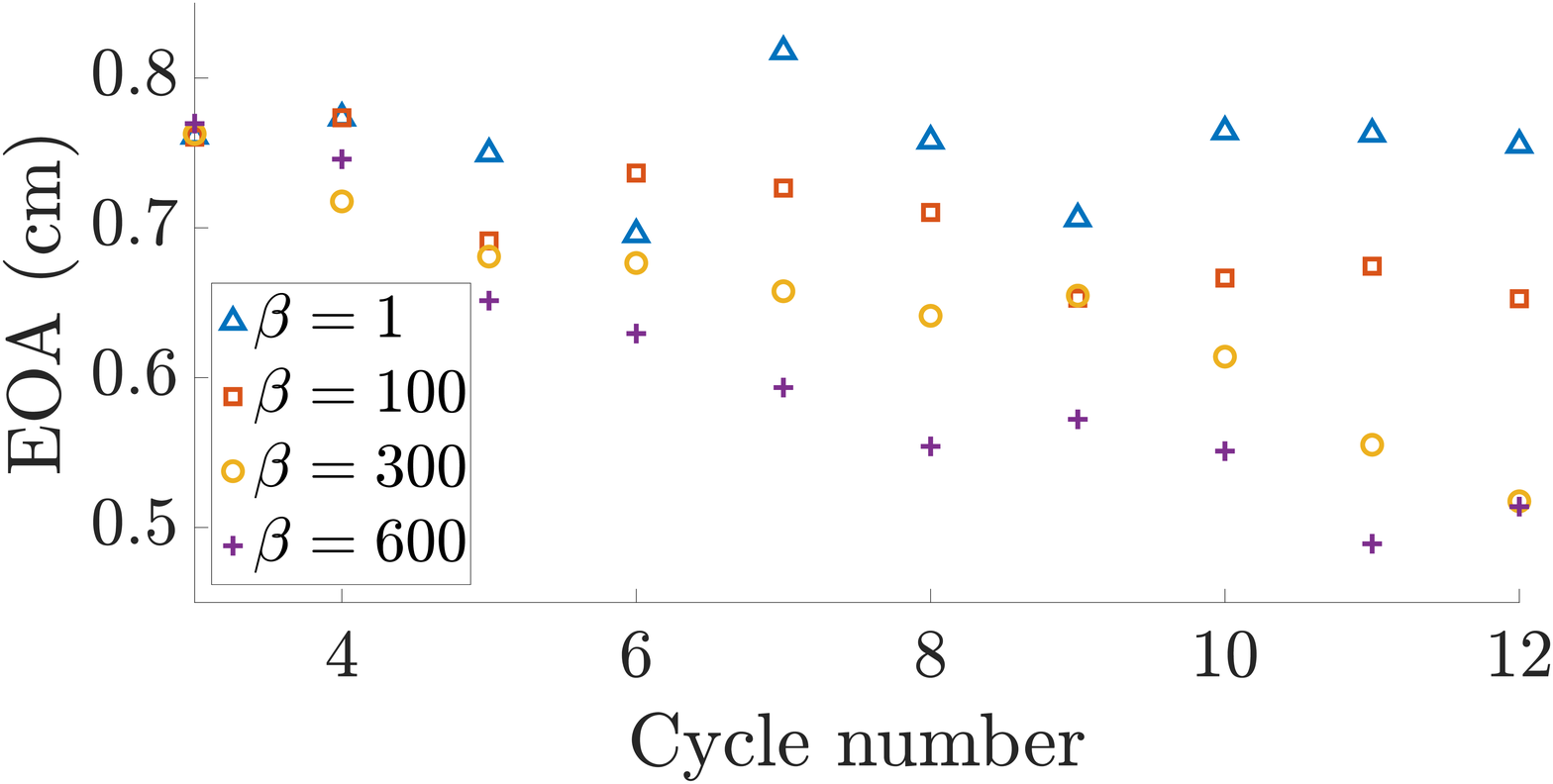}
\subfigimg[width=0.5\columnwidth,hsep=-0.75em,vsep=2em,pos=ul]{(b)}{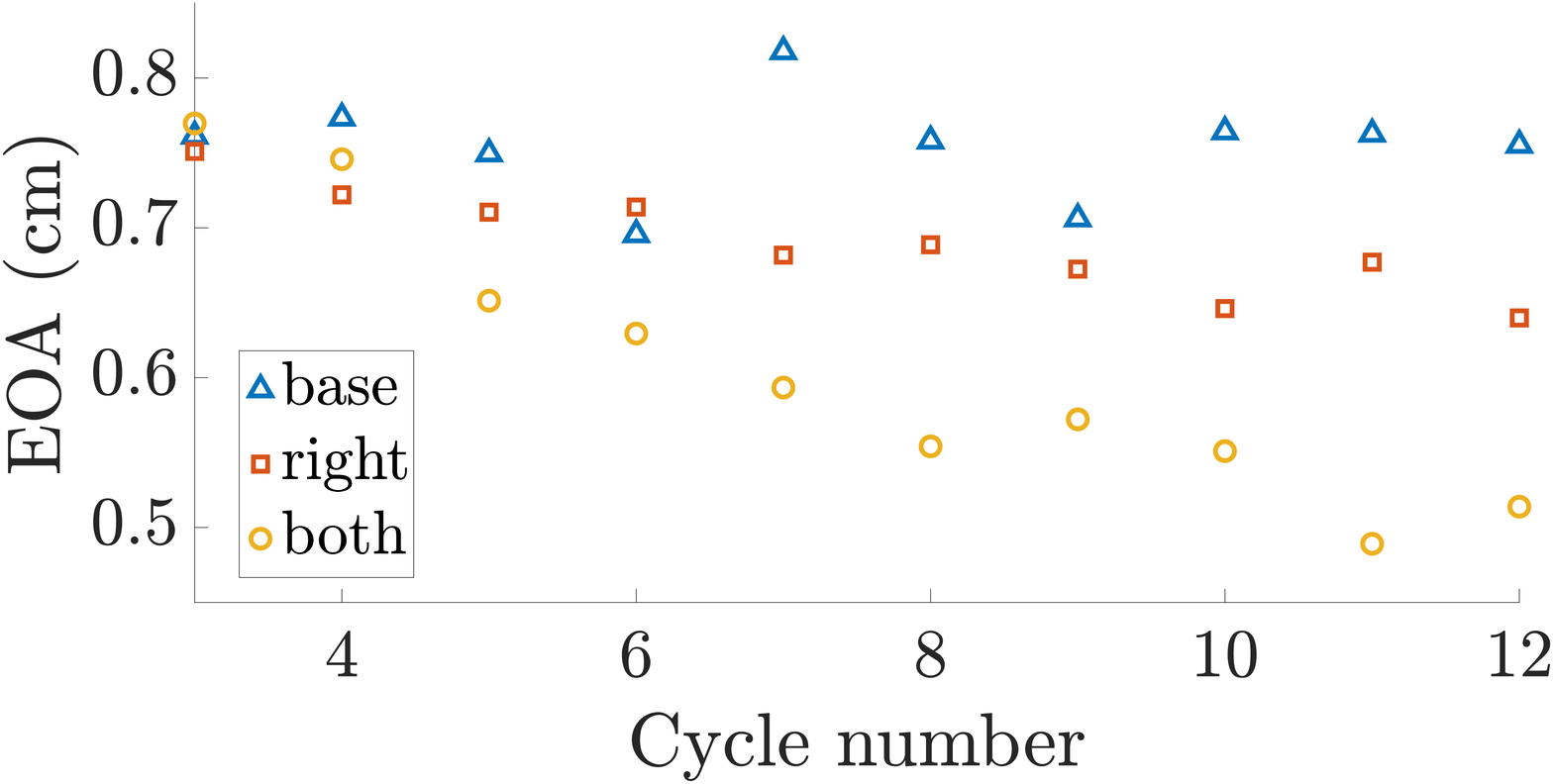}
\caption{The effective orifice area (EOA) for each cycle. Panel (a) shows the EOA with $\beta = 1$, 100, 300, and 600 for accumulation on both leaflets. Note that $\beta = 1$ corresponds to a baseline where the leaflet stiffness is constant. Panel (b) shows the EOA for accumulation on no leaflets (base), accumulation on the right leaflet (right), and accumulation on both leaflets (both). We observe a general decrease of EOA over time as the accumulation increases. There is a greater decrease of EOA if accumulation occurs on both leaflets than just the right leaflet.}
\label{fig:eoa}
\end{figure}

We additionally compute the effective orifice area (EOA). The EOA $A_\text{AV}$ is computed using conservation of mass in \cref{eq:nes_1} by assuming the relation $V_\text{LVOT}A_\text{LVOT} = V_\text{AV}A_\text{AV},$ in which $V_\text{LVOT}$ is the average time integral flow rate through the left ventricle outflow tract during each cycle when the valve is open, $A_\text{LVOT}$ is the area of the left ventricle outflow tract, and $V_\text{AV}$ is the average time integral flow rate through the aortic valve. $V_\text{LVOT}$ is computed from the boundary condition model described in \cref{sec:bdry_conds}. To compute $V_\text{AV}$, we interpolate the velocity to the midpoint between the two points on each leaflet that are closest during systole. $V_\text{AV}$ is then computed as the time integral of the component of this interpolated velocity normal to the valve ring. \Cref{fig:eoa} plots the EOA for each cycle. We observe a general decrease in EOA as more the total surface concentration $\Cb$ increases. This indicates that the fluid velocity increases to compensate for the stiffening of the valve. The EOA decreases more when accumulation occurs on both leaflets when compared to accumulation on only the right leaflet. Because the left leaflet remains at a constant stiffness, the left leaflet opens more to compensate for the stiffening of the right leaflet. This causes the jet to shift towards the left leaflet, as shown in \cref{fig:during_sim_right}.

\section{Conclusions}
This study presents new numerical methods incorporating both deposition and fluid-structure interaction to simulate leaflet thrombosis. The simplified thrombosis model serves as a stepping stone to demonstrate the capabilities of our simulation approach that includes concentration fields describing fluid-phase platelets and structure-bound platelets. Platelets can deposit onto the leaflet surface, and bound platelets can dissociate into the fluid. In our model, the stiffness of the leaflet is a function of the bound platelet concentration. We have shown that our model is capable of realizing drops in pressure and decreases in effective orifice area, without fully occluding the aortic valve. The results also show that the stiffness of the valve can lead to a variety of flow features in the sinus of Valsalva region. These flow features affect the amount of material that is locally present to deposit over the leaflets.

Extensions of this model to three dimensions require an efficient method for solving the advection-diffusion equation in complex, time-evolving domains. The method utilized here requires the computation of cut-cell volumes and intersections, which remain challenging in three spatial dimensions. Recent approaches to this class of problems include mesh-free RBF-FD methods \cite{Shankar2018} and volume penalization methods \cite{Thirumalaisamy2022}. The implementation of a more physiological model of thrombosis remains important future work. A primary roadblock is the disparate time scales present in thrombosis. While the heart beats on the order of seconds, blood clots can form in hours to days. The use of conditionally stable time stepping limits the numerical methods to time steps that resolve the fastest timescale, which in this model is that of the fluid-structure interaction. Recent work in multiscale time stepping algorithms \cite{Frei2020,Frei2016} could enable extensions of our modeling framework to enable such long-time simulations. Further, with multiscale time stepping algorithms, this model could be extended to study the affect of saturation of the bound concentration field. While platelet deposition is important and the beginning step of thrombus formation, a significant portion of the clot may be from coagulation and fibrin mesh formation. However, a complete model of thrombosis will require a computational model in which the blood clot on the moving valve leaflets grows into the fluid \cite{Fogelson2008,Leiderman2011}. The development of such a model that incorporates FSI is ongoing.

The new approaches described herein should be considered an important steppingstone for thrombosis models in many different contexts. This model is the first of its kind to incorporate both adhesion of a surface concentration to the surface of the leaflets and feedback into the fluid-structure interaction. Further, this model can be adapted to model deposition or absorption along other moving boundaries, such as for particulate flow in the lungs or drug absorption in the gut.

\section*{Acknowledgements}
We acknowledge funding from the NIH (Awards U01HL143336 and R01HL157631) and NSF (OAC 1652541 and OAC 1931516).

\printbibliography[heading=bibintoc]

\end{document}